\documentclass[11pt,a4paper]{article}
\usepackage{epsfig}
\usepackage{graphics}
\usepackage{amsmath,amssymb}
\usepackage{subfigure}
\usepackage[T1]{fontenc}
\usepackage{graphicx}
\usepackage{epstopdf}
\usepackage{color}

\newcommand{\thefont}[2]{\fontsize{#1}{#2}\fontshape{n}\selectfont}
\newcommand{\1}{\rlap{\thefont{10pt}{12pt}1}\kern.16em\rlap{\thefont{11pt}{13.2pt}1}\kern.4em}
\newcommand{\RR}{\ensuremath{\mathbb{R}} }
\newcommand{\CC}{\ensuremath{\mathbb{C}} }
\newcommand{\ZZ}{\ensuremath{\mathbb{Z}} }
\newcommand{\NN}{\ensuremath{\mathbb{N}} }

\newcommand{\EE}{\ensuremath{{\mathbb E}}}
\newcommand{\opO}{\ensuremath{\mathcal O}}

\def\E{\mathbb{E}_{\theta}}

\newtheorem{ass}{Assumption}
\newtheorem{theo}{Theorem}

\newtheorem{lemma}{Lemma}

\newtheorem{coro}{Corollary}

\numberwithin{equation}{section} \numberwithin{ass}{section}
\numberwithin{theo}{section} \numberwithin{prop}{section}
\numberwithin{lemma}{section} \numberwithin{definition}{section}
\numberwithin{rmq}{section}

\title{Sharp template estimation in a shifted curves model}

\author{{\em J\'er\'emie Bigot, S\'ebastien Gadat  \&  Cl\'ement Marteau} \\
Institut de Math\'ematiques de Toulouse\\
Universit\'e de Toulouse et CNRS (UMR 5219)\\
31062 Toulouse, Cedex 9, France\\
{\small {\tt \{Jeremie.Bigot,Sebastien.Gadat,Clement.Marteau\}@math.univ-toulouse.fr, } }}

\begin{document}

\maketitle

\thispagestyle{empty}

\begin{abstract}
This paper considers the problem of adaptive estimation of a template in a randomly shifted curve model. Using the Fourier transform of the data, we show that this problem can be transformed into a stochastic linear inverse problem.  Our aim is to approach the estimator that has the smallest risk on the true template over a finite set of linear estimators defined in the Fourier domain. Based on the principle of unbiased empirical risk minimization, we derive a nonasymptotic oracle inequality in the case where the law of the random shifts is known. This inequality can then be used to obtain adaptive results on Sobolev spaces as the number of observed curves tend to infinity. Some numerical experiments are given to illustrate the performances of our approach.
\end{abstract}


\section{Introduction}
\subsection{Model and objectives}
The goal of this paper is to study a special class of stochastic inverse problems. We consider the problem of estimating a curve $f$, called template or shape function, from the observations of $n$ noisy and randomly shifted curves  $Y_1, \dots Y_n$ coming from the following Gaussian white noise model:
\begin{equation}\label{model}
dY_j(x) = f(x-\tau_j) dx + \epsilon dW_j(x), \; x \in [0,1], \; j=1,\ldots,n
\end{equation}
where $W_{j}$ are independent standard Brownian motions on $[0,1]$, $\epsilon$ represents a level of noise common to all curves, the $\tau_j$'s are unknown random shifts, $f$ is the unknown template to recover, and $n$ is the number of observed curves that may be let going to infinity to study asymptotic properties. This model is realistic in many situations where it is reasonable to assume that the observed curves represent replications of almost the same process and when a large source of variation in the experiments is due to transformations of the time axis. Such a model is commonly used in many applied areas dealing with functional data such as neuroscience (see e.g. \cite{IRT}) or biology (see e.g.  \cite{ronn}). A well known problem in functional data analysis is the alignment of similar curves that differ by a time transformation to extract their common features, and (\ref{model}) is a simple model where $f$ represents such common features (see \cite{ramsil}, \cite{ramsil2}  for a detailed introduction to curve alignment problems in statistics).

The function $f : \RR \to \RR$ is assumed to be of period $1$ so that the model (\ref{model}) is well defined, and the shifts  $\tau_j$ are supposed to be  independent and identically distributed (i.i.d.) random variables with density $g : \RR \to \RR$ with respect to the Lebesgue measure $dx$ on $\RR$. Estimating $f$ can be seen as a stochastic inverse problem as this template is not observed directly, but through $n$ independent realizations of the stochastic operator $A_{\tau} : L^{2}_{p}([0,1]) \to L^{2}_{p}([0,1])$ defined by
$$
A_{\tau}(f)(x) = f(x-\tau), \; x \in [0,1],
$$
where $L^{2}_{p}([0,1])$ denotes the space of squared integrable functions on $[0,1]$ with period 1, and $\tau$ is  random variable with density $g$. The additive Gaussian noise makes this problem ill-posed, and \cite{BG} have shown that estimating $f$ in such models is in fact a deconvolution problem where the density $g$ of the random shifts plays the role of the convolution operator. For the $L^{2}$ risk on $[0,1]$, \cite{BG} have derived the minimax rate of convergence for the estimation of $f$ over Besov balls as $n$ tends to infinfity. This minimax rate depends both on the smoothness of the template and on the decay of the Fourier coefficients of the density $g$. This is a well known fact for standard deterministic deconvolution problem in statistics, see e.g. \cite{fan}, \cite{donoho}, but the results in \cite{BG} represent a novel contribution and a new point of view on template estimation in stochastic inverse problems such as (\ref{model}).

However, the approach followed in \cite{BG}  is only asymptotic, and the main
goal of this paper is to derive non-asymptotic results  to study the estimation
of $f$  by keeping fixed the number $n$ of observed curves.

\subsubsection{Deconvolution formulation}

Let us first explain how the model  (\ref{model}) can be transformed into a
deconvolution problem as the one studied in \cite{DJKP95jrssb}. Denote $G$ the
following density function defined on $[0;1]$ as
$$
G(x) = \sum_{k \in \mathbb{Z}} g(x+k).
$$
The density $G$ exists as soon as $g$ satisfies the weak condition $g(x) \leq
\frac{C}{1+|x|^{\nu}}$ for any $\nu>1$ and suitable constant $C$. Note that the
Fourier coefficients of G are given by
$$
\int_{0}^1 G(t) e^{- i 2 \pi l t} dt = \int_{- \infty}^{\infty} g(t) e^{-i 2
\pi l t} dt = \gamma_l
$$
Consider now the 1-periodization of $f$ extended to $\mathbb{R}$, one has
$$
\int_{0}^1 f(x-\tau) G(\tau) d\tau = \int_{- \infty}^{\infty} f(x-\tau) g(\tau)
d\tau.
$$
The observations $Y_j$ can be written as
\begin{equation}\label{eq:model-inverse}
dY_j(x) = f \star G(x) dx + \xi_j(x) dx + \epsilon dW_j(x), \end{equation}
where $\xi_j$ is a second noise term defined as $ \xi_j(x) = f(x-\tau_j) - f
\star G(x)$. Hence, our model can be seen as a deconvolution problem with a
noisy operator $H: f \mapsto f \star G+ \xi$ and a more classical independent
additive noise $W$. Note also that the realizations $H_j:f \mapsto f \star G+
\xi_j $ are unbiased realizations of the operator $H$ but presents a variance
term which depends on the function $f$ we want to estimate. This appears to be
a new setting in the field of inverse problem with unknown operators as considered in \cite{cavheng}, \cite{EfromovichK01}, \cite{hr}, \cite{marteau1} and \cite{cavraim}.

We will see in the sequel that the additive noise $\xi$ which depends on $f$
slightly modifies the quadratic risk and the way to estimate $f$ when compared to classical procedures used in standard inverse problems with a deterministic operator.

\subsection{Fourier Analysis and an inverse problem formulation}

 Supposing that $f \in L_{p}^{2}([0,1])$, we denote by $\theta_k$ its $k^{th}$ Fourier coefficient, namely:
$$
\theta_k = \int_{0}^1 e^{- 2 i k \pi x} f(x) dx.
$$
In the Fourier domain, the model (\ref{model}) can be rewritten as
\begin{equation} \label{ckj}
c_{j,k} :=  \int_{0}^1 e^{- 2 i k \pi x} dY_{j}(x) =  \theta_k e^{-i 2 \pi k \tau_j} + \epsilon  z_{k,j}
\end{equation}
where $z_{k,j}$ are i.i.d.  $\mathcal{N}_{\CC}\left(0,1\right)$ variables, i.e. complex Gaussian variables with zero mean and such that $\EE |z_{k,j}|^{2} = 1$. This means that the real and imaginary parts of the  $z_{k,j}$ 's are Gaussian variables with zero mean and variance 1/2. Thus, we can compute the sample mean of the $k^{th}$ Fourier coefficient over the $n$ curves as
\begin{equation} \label{eqseqmodel}
\tilde{c}_{k}:= \frac{1}{n} \sum_{j=1}^n c_{k,j} = \theta_k \tilde{\gamma}_{k} + \frac{\epsilon}{\sqrt{n}} \xi_k,
\end{equation}
where
\begin{equation} \label{eq:gammatilde}
\tilde{\gamma}_{k} := \frac{1}{n} \sum_{j=1}^n e^{-i 2 \pi k \tau_j},
\end{equation}
and the $\xi_k$'s are i.i.d. complex Gaussian variables with zero mean and variance $1$.  The Fourier coefficients $\tilde c_k$ in equation (\ref{eqseqmodel}) can be viewed as observations  coming from a statistical inverse problem. Indeed, the standard sequence space model of an ill-posed statistical inverse problem is (see \cite{cavgopitsy} and the references therein)
\begin{equation} \label{eqinvpb}
c_{k} = \theta_{k} \gamma_{k} + \sigma z_{k},
\end{equation}
where the $\gamma_{k}$'s are eigenvalues of a known linear operator, $z_{k}$ are random noise variables and $\sigma$ is a level of noise which goes to zero for studying asymptotic properties. The issue in such models is to recover the coefficients $\theta_{k}$ from the observations $c_{k}$ under various conditions on the decay to zero of the  $\gamma_{k}$'s as $|k| \to + \infty$. A large class of estimators for the problem (\ref{eqinvpb}) can be written as
$$
\hat{\theta}_{k} = \lambda_{k}  \frac{c_{k}}{\gamma_{k}},
$$
where $\lambda = (\lambda_{k})_{k \in \ZZ}$ is a sequence of reals called filter. Various estimators of this form have been studied in a number of papers, and we refer to \cite{cavgopitsy} for more details.

In a sense, we can view equation (\ref{eqseqmodel}) as an inverse problem (with $\sigma = \frac{\epsilon}{\sqrt{n}}$) where the eigenvalues of the linear operator are the Fourier coefficients of the density $g$ of the shifts i.e.
$$
\gamma_k := \mathbb{E} \left( e^{- i 2 \pi k \tau} \right) = \int_{-\infty}^{+\infty} e^{- i 2 \pi k x} g(x) dx .
$$

Indeed, let us assume that the density $g$ of the random shifts is known. In this case, to estimate the Fourier coefficients of $f$, one can perform a deconvolution step of the form
\begin{equation}
\hat{\theta}_{k} = \lambda_{k} \frac{\tilde{c}_{k}}{\gamma_{k}} \label{eq:deftheta},
\end{equation}
where $\tilde c_k$ is defined in (\ref{eqseqmodel}) and  $\lambda = (\lambda_{k})_{k \in \ZZ}$ is a filter whose choice will be discussed later on.
Theoretical properties and optimal choices for the filter $\lambda$ will be presented in the case where the coefficients $\gamma_{k}$ are known. Such a framework is commonly used  in inverse problems such as (\ref{eqinvpb}) to obtain consistency results and to study asymptotic rates of convergence, where it is generally supposed that the law of the additive error is Gaussian with zero mean and {\it known} variance $\sigma^{2}$, see e.g \cite{cavgopitsy}. In model (\ref{model}), the random shifts may be viewed as a second source of noise and for the theoretical analysis of this problem the law of this other random noise is also supposed to be known.

Recently, some papers have addressed the problem of regularization with partially known operator. For instance, \cite{cavheng} consider the case where the eigenvalues are unknown but independently observed. They deal with the model:
\begin{equation}
c_k = \gamma_k \theta_k + \epsilon \xi_k, \ \tilde \gamma_k = \gamma_k + \sigma \eta_k, \ \forall k\in \NN,
\label{eq:modelLN}
\end{equation}
where $(\xi_k)_{k\in\NN}$ and $(\eta_k)_{k\in\NN}$ denote i.i.d standard gaussian variables. In this case, each coefficient $\theta_k$ can be estimated by $\tilde \gamma_k^{-1} c_k$. Similar models have been considered in \cite{cavraim}, \cite{marteau1} or \cite{marteau2}. In a more general setting, we may refer to \cite{EfromovichK01} and \cite{hr}.

In this paper, our framework is sligthly different in the sense that the operator is stochastic, but the regularization is operated using deterministic
eigenvalues. Hence the approach followed in the previous papers is no directly
applicable to model (\ref{model}). We believe that estimating $f$ in model
(\ref{model}) without the knowledge of $g$ remains a difficult task, and this
paper is a first step to address this issue.

\subsection{Previous work in template estimation and shift recovery}

The problem of estimating the common shape of a set of curves that differ by a time transformation is usually referred to as the curve registration problem, and it has received a lot of attention in the literature over the last two decades.  Among the various methods that have been proposed, one can distinguish between landmark-based approaches which aim at aligning common structural points of the curves (typically locations of extrema) see e.g \cite{gaskneip}, \cite{kneipgas}, \cite{big}, and nonparametric modeling  of the warping functions to align a set of curves see e.g \cite{ramli}, \cite{wanggas}, \cite{liumuller}. However, in these papers, studying consistent estimates of the common shape $f$ as the number of curves $n$ tends to infinity is generally not considered.

In the simplest case of  shifted curves, various approaches have been developed. Self-modelling regression methods proposed by \cite{kg} are semiparametric models where each observed curve is a parametric transformation of a common regression function. Such models are usually referred to as shape invariant models and estimation in this setting is usually done by iterating the following two steps: estimation of the parameters of the transformations (here the shifts) given a reference curve, and nonparametric estimation of a template by aligning the observed curves given a set of known transformation parameters.  \cite{kg}  studied the consistency of such a two steps procedure in an asymptotic framework where both the number of functions $n$ and the number of observed points per curves grows to infinity. Due to the asymptotic equivalence between the white noise model and nonparametric regression with an equi-spaced design (see \cite{BM96aos}), such an asymptotic framework in our setting would correspond to the case where both $n$ tends to infinity and $\epsilon$ is let going to zero. In this paper we prefer to focus only on the case where  $n$ may be let going to infinity, and to leave fixed the level of additive noise in each observed curve.

Based on a model with curves observed at discrete time points, semiparametric estimation of the shifts and the shape function is proposed in \cite{mazaloubgam} and \cite{vimond} as the number of observations per curve grows, but with a fixed number $n$ of curves. A generalisation of this approach for the estimation of scaling, rotation and translation parameters for two-dimensional images is also proposed in \cite{BGV}, but also with a fixed number of observed images.  Semiparametric and adaptive estimation of a shift parameter in the case of a single observed curve in a white noise model is also considered by \cite{DGT} and \cite{D}. Estimation of a common shape for randomly shifted curves and asymptotic in $n$ is considered in \cite{ronn} from the point of view of semiparametric estimation when the parameter of interest is infinite dimensional.

However, in all the above cited papers rates of convergence or oracle inequalities for the estimation of the template are generally not studied. Moreover, our procedure differs from the approaches classically used in curve registration as our estimator is obtained in only one very simple step, and it is not based on an alternative scheme between estimation of  the shifts and averaging of back-transformed curves given estimated values of the shifts parameters.

Finally, note that \cite{loubcast} and \cite{IRT} consider a model similar to
(\ref{model}), but they rather  focus on the the estimation of the density $g$
of the shifts as $n$ tends to infinity. Using such an approach could be a good
start for studying the estimation of the template $f$ without the knowledge of
$g$. However, we believe that this is far beyond the scope of this paper, and
we prefer to leave this problem open for future work.

\subsection{Organization of the paper}

In Section \ref{sec:estimf}, we consider an estimator of the shape function $f$ based on spectral cut-off when the eigenvalues $\gamma_{k}$ are known. Based on the principle of unbiased risk minimization developed by \cite{cavgopitsy}, we derive an oracle inequality that is then used to derive an adaptive estimator of $f$ on Sobolev spaces. This estimator is based on the Fourier transform of the curves with a data-based choice of the frequency cut-off. In Section \ref{sec:minimax}, we study asymptotic properties of this estimator in terms of minimax rates of converge over Sobolev balls. Finally in Section \ref{sec:simu}, a short simulation study is proposed  to illustrate the numerical properties  of the estimator. All proofs are deferred to a technical section at the end of the paper.

\section{Estimation of the common shape} \label{sec:estimf}

In the following, we assume that the Fourier coefficients  $\gamma_k$ are known. In this situation it is possible to choose a data-dependent filter $\lambda^{\star}$ that mimic the performances of an optimal filter  $\lambda^{0}$ called oracle that would be obtained if we knew the true template $f$. The performances of this filter are related to the performances of the filter  $\lambda^{0}$ via an oracle inequality. In this section, most of our results are non-asymptotic and are thus related to the approach proposed in \cite{cavgopitsy} to study standard statistical inverse problems via oracle inequalities.

\subsection{Smoothness assumptions for the density $g$}

In a deconvolution problem, it is well known that  the difficulty of estimating
$f$ is quantified by the decay to zero of the $\gamma_{k}$'s as $|k| \to +
\infty$. Depending how fast these Fourier coefficients tend to zero as $|k| \to
+ \infty$, the reconstruction of $f$ will be more or less accurate. This
phenomenon was systematically studied by \cite{fan} in the context of density
deconvolution. In this paper, the following type of assumption on  $g$
is considered:

\begin{ass}\label{assordi}
The Fourier coefficients of $g$ have a polynomial decay i.e.  for some real $\beta \geq 0$, there exists two constants $C_{max} \geq C_{min} > 0$ such that for all $k \in \ZZ$
\begin{equation} \label{eq:polydecay}
C_{min} |k|^{-\beta} \leq  |\gamma_{k}| \leq C_{max} |k|^{-\beta}.
\end{equation}
\end{ass}
Remark that the knowledge of the constants $C_{max}, C_{min}$ and $\beta$ will
not be necessary for the construction of our estimator.

\subsection{Risk decomposition}

Assuming that  $\gamma_k \neq 0$ for all $k \in \ZZ$,  we recall that an
estimator of the $\theta_{k}$'s is given by, see equation (\ref{eq:deftheta})
$$
\hat{\theta}_{k} = \lambda_{k} \frac{\tilde{c}_{k}}{\gamma_{k}}
$$
where $\lambda = (\lambda_{k})_{k \in \ZZ}$ is a real sequence. Examples of
commonly used filters include projection weights $\lambda_{k} = \1_{|k| \leq
N}$ for some integer $N$, and the Tikhonov weights $\lambda_{k} =
1/(1+(|k|/\nu_{2})^{\nu_{1}})$ for some parameters $\nu_{1} > 0$ and $\nu_{2} >
0$. Based on the $\hat{\theta}_k$'s, one can estimate the signal $f$ using the
Fourier reconstruction formula.

The problem is then to choose the sequence $(\lambda_{k})_{k \in \ZZ}$ in an optimal way with respect to an appropriate risk. For a given filter $\lambda$  we use the classical $\ell_{2}$-norm to define the risk of the estimator $\hat{\theta}(\lambda) = (\hat{\theta}_k)_{k \in \ZZ}$
\begin{equation} \label{eq:defrisk}
R(\theta,\lambda) = \mathbb{E} \| \hat{\theta}(\lambda) - \theta\|^2 = \mathbb{E} \sum_{k \in \mathbb{Z}} |\hat{\theta}_k -  \theta_k|^{2}
\end{equation}

Note that analyzing the above risk (\ref{eq:defrisk}) is equivalent to analyze the mean integrated square risk $R(\hat{f}_{\lambda},f) = \EE \|\hat{f}_{\lambda} - f \|^{2} =  \EE \left(  \int_{0}^{1} (\hat{f}_{\lambda}(x)-f(x) )^{2} dx \right)$ for the estimator $\hat{f}_{\lambda}(x) = \sum_{k \in \ZZ} \hat{\theta}_k e^{- 2 i k \pi x}$. The following lemma gives the bias-variance decomposition of
$R(\lambda,\theta)$.

\begin{lemma} \label{lem:decomp}
For any given nonrandom filter $\lambda$,  the risk of the estimator $\hat{\theta}(\lambda)$ can be decomposed as
\begin{equation} \label{eq:decomprisk}
R(\theta, \lambda) = \underbrace{\displaystyle\sum_{k \in \mathbb{Z}}
(\lambda_k-1)^2 |\theta_k|^2}_{Bias} + \underbrace{ \frac{1}{n}
\displaystyle\sum_{k \in \mathbb{Z}}
\lambda_k^2\frac{\epsilon^2}{|\gamma_k|^2}}_{V_1} +
\underbrace{\frac{1}{n}\sum_{k \in \mathbb{Z}}\left[\lambda_k^2|\theta_k|^2
\left(\frac{1}{|\gamma_k|^2} -1 \right) \right]}_{V_2}
\end{equation}
\end{lemma}
For a fixed number of curves $n$ and a given shape function $f$, the problem of
choosing an optimal filter in a set of possible candidates  is to find the best
tradeoff between low bias and low variance in the above expression. However,
this decomposition does not correspond exactly to the classical bias-variance
decomposition for linear inverse problems. Indeed, the variance term in
(\ref{eq:decomprisk}) is the sum of two terms and differs from the classical
expression of the variance for linear estimator in statistical inverse
problems. Using our notations, the classical variance term is
$V_1=\frac{\epsilon^2}{n} \displaystyle\sum_{k \in \mathbb{Z}}
\frac{\lambda_k^2}{|\gamma_k|^2} $ and appears in most of linear inverse
problems.

However, contrary to standard inverse problems, the variance term of the risk
also depends on the Fourier coefficients $\theta_k$ of the unknown function $f$
to recover. Indeed, our data $ \gamma_k ^{-1} \tilde c_k$ are noisy
observations of $\theta_k$:
$$
\gamma_k ^{-1} \tilde c_k = \theta_k +
\left(\frac{\tilde\gamma_k}{\gamma_k} -1 \right) \theta_k + \frac{\epsilon}{\sqrt{n}} \gamma_k^{-1} \xi_k,
$$
and we invert the problem using the sequence $(\gamma_k)_{k\in\NN}$ instead of
$(\tilde \gamma_k)_{k\in\NN}$, which is involved in the construction of the coefficient
$c_k$.
It explains the presence of the second term $V_2$. In particular, the quadratic risk is expressed in its usual form in the case where $\tilde\gamma_k=\gamma_k$. \\

A similar phenomenon occurs with the model (\ref{eq:modelLN}), although it is more difficult to quantify. Indeed, in this setting:
$$ \tilde \gamma_k^{-1} c_k = \theta_k + \left( \frac{\gamma_k}{\tilde\gamma_k} -1 \right) \theta_k + \epsilon \tilde\gamma_k^{-1} \xi_k, \ \forall k\in\NN.$$
Hence, we also observe an additionnal term depending on $\theta$. This term is controled using a Taylor expension but the quadratic risk cannot be expressed in a simple form. We refer to \cite{marteau2} for a discussion with some numerical simulation and to \cite{cavheng}, \cite{EfromovichK01}, \cite{hr}, \cite{marteau1} and \cite{cavraim}.

\subsection{An oracle estimator and unbiased estimation of the risk}

Suppose that one is given a finite set of possible candidate filters $\Lambda = (\lambda^{N})_{N\in \ I}$, with $\lambda^{N}=(\lambda^{N}_{k})_{k \in \ZZ}, \; N\in I \subset \NN$ which satisfy some general conditions to be discussed later on. In the case of projection filters, $\Lambda$ can be for example the set of filters $\lambda^{N}_{k} = \1_{|k| \leq N}, k \in \ZZ$ for $N = 1,\ldots,m_{0}$. Given a set of filters $\Lambda$, the best estimator corresponds to the filter $\lambda^0$, called oracle, which minimizes the risk $R(\lambda,\theta)$ over $\Lambda$ i.e.
\begin{equation}
\label{eq.lambda0} \lambda^0 := \arg \min_{\lambda \in \Lambda}
R(\lambda,\theta).
\end{equation}
This filter is called an oracle because it cannot be computed in practice as the sequence of coefficients $\theta$ is unknown.
However, the oracle $\lambda^0$ can be used as a benchmark to evaluate the quality of a data-dependent filter $\lambda^{\star}$ chosen in the set $\Lambda$. This is the main interpretation of the oracle inequality that we will develop in the next section.

Now, suppose that it is possible to construct an unbiased estimator $\hat{\Theta}^2_k$ of $|\theta_k|^2$. For any nonrandom filter $\lambda$, using $\hat{\Theta}^2_k$, one can compute an estimator $\tilde{U}(\lambda,X)$ of the risk $R(\lambda,\theta)$.
Then, for choosing a data-dependent filter, the principle of unbiased risk estimation (see \cite{cavgopitsy} for further details) simply suggests to minimize the criterion $U(\lambda,X)$ over $\lambda \in \Lambda$ instead of the criterion $R(\lambda,\theta)$.  Our data-dependent choice of $\lambda$ is thus
\begin{equation}
\label{eq.lambda} \lambda^{\star} := \arg \min_{\lambda \in \Lambda} \tilde{U}(\lambda,X).
\end{equation}
Typically, in practice, all the filters $\lambda \in \Lambda$ are such that
$\lambda_{k}= 0$ (or vanishingly small) for all $k$ large enough.  Hence, for
such choices of filters, numerical computation of the above expression is thus
feasible since it only involves the computation of finite sums.

\subsection{Oracle inequalities for projection filters}

\subsubsection{Unbiased Risk Estimation (URE)}
For the sake of simplicity, we only consider spectral cut-off schemes in the following. In this case, $\Lambda$ corresponds to the set of filters $(\1_{|k| \leq N})_{k \in \ZZ}$ for $N \in \NN$. All the results presented in this paper could be generalized to wider families of estimators (Tikhonov, Landweber, Pinsker,...). The price to pay is to get longer and more technical proofs.

From Lemma \ref{lem:decomp}, the quadratic risk $R(\theta,\lambda):= R(\theta,N)$ of a projection filter can be written as:
\begin{eqnarray*}
R(\theta,N) &=& \sum_{|k| > N}  |\theta_k|^2 +   \displaystyle
\frac{\epsilon^2}{n}  \sum_{|k| \leq N} |\gamma_k|^{-2} +  \frac{1}{n}
\sum_{|k| \leq N}  |\theta_k|^2 \left(\frac{1}{|\gamma_k|^2} -1 \right) \\
& = & \|\theta\|_2^2 - \sum_{|k| \leq N}  |\theta_k|^2 + \displaystyle
\frac{\epsilon^2}{n}  \sum_{|k| \leq N} |\gamma_k|^{-2} +  \frac{1}{n}
\sum_{|k| \leq N}  |\theta_k|^2 \left(\frac{1}{|\gamma_k|^2} -1 \right)
\end{eqnarray*}

We aim to minimize $R$ with respect to $N$ while $\theta$ is unknown.
 Using $\hat{\Theta}^2_k = \gamma_k^{-2} \left[ |\tilde c_k|^2 -
\frac{\epsilon^2}{n}\right]$ as an unbiased estimator of $|\theta_k |^2$, we
minimize U defined as
\begin{equation}
U(Y,N)= - \left(1-\frac{1}{n}\right)\sum_{|k|\leq N} |\gamma_k|^{-2} \left\lbrace
|\tilde c_k|^2 - \frac{\epsilon^2}{n} \right\rbrace + \frac{\epsilon^2}{n}
\sum_{|k|\leq N} |\gamma_k|^{-2} + \frac{1}{n} \sum_{|k|\leq N} |\gamma_k
|^{-4} | \left\lbrace |\tilde c_k|^2 -\frac{\epsilon^2}{n} \right\rbrace,
\label{def:urisque1}
\end{equation}
which is an unbiased risk estimator of $R(\theta,N) - \|\theta\|_2^2$.

 Unfortunately, such a criterion does not lead to satisfying results.
Instead of the approach developed in \cite{cavheng}, we take into account the
error generated by the use of an approximation of the eigenvalues. The
estimator related to the criterion (\ref{def:urisque1}) involves processes that
require a specific treatment. In order to contain these processes, we will
consider in the following the criterion
\begin{equation}
\bar U(Y,N)= - \sum_{|k|\leq N} |\gamma_k|^{-2} \left\lbrace |\tilde c_k|^2 -
\frac{\epsilon^2}{n} \right\rbrace + \frac{\epsilon^2}{n} \sum_{|k|\leq N}
|\gamma_k|^{-2} + \frac{\log^2(n)}{n} \sum_{|k|\leq N} |\gamma_k |^{-2} 
\left\lbrace |\bar c_k| -\frac{\epsilon^2}{n} \right\rbrace, \label{def:u1pen}
\end{equation}
Remark that $\bar U(Y,N)$ can be written as $U(Y,N)+ \mathrm{pen}(N)$ where
$(\mathrm{pen}(N))_{N\in \NN}$ denotes a penalty term. It appears from the
proofs that this penalty is a natural candidate for the control of the
processes involved in the behavior of the estimator constructed below. The
associated data-based filter is defined as
\begin{equation}
N^{\star} = \mathrm{arg} \min_{N\leq m_0} \bar U(Y,N),
\label{def:fenetre}
\end{equation}
where
\begin{equation}
m_0 = \inf \left\lbrace k: |\gamma_k|^2 \leq \frac{\log ^2 n}{n} \right\rbrace
- 1. \label{def:m0}
\end{equation}
Remark that we do not minimize our criterion $\bar U(Y,N)$ over $\NN$ but rather for $N\leq
m_0$. Indeed, each coefficient $\theta_k$ is estimated by $\gamma_k^{-1} \tilde c_k$ where
$\gamma_k = \EE [\tilde \gamma_k]$. Hence, the ratio $\gamma_k^{-1}\tilde \gamma_k$ should be as close as possible to 1. Since $\gamma_k\rightarrow 0$ as $k\rightarrow +\infty$ and the variance of $\tilde \gamma_k$ is constant in $k$, it seems clear that large $k$ should be avoided.

Similar bounds on the resolution level are used in papers related to partially known operator: see for instance
\cite{cavheng} or \cite{EfromovichK01}. This bounds have to be carefully chosen but are not of first importance. In general,
estimating the operator is easier than estimating the function $f$.

\subsubsection{Sharp estimator of the risk}

We are now able to propose a first adaptive estimator. In the following, we denote by $\theta^{\star}$ the estimator related to the bandwidth $N^{\star}$ namely
\begin{equation}
\theta^{\star}_{k} = \frac{\tilde{c}_{k}}{\gamma_{k}} \1_{\lbrace k \leq N^{\star}\rbrace} . \label{def:thetastar} 
\end{equation}
The next theorem summarizes the performances of $\theta^{\star}$ through a simple oracle inequality. The proof is postponed to the Section 5.

\begin{theo}
\label{th:full}
Let $\theta^{\star}$ defined by (\ref{def:thetastar}) and assume that the density $g$ satisfies Assumption \ref{assordi}. Then, there exists $0<\gamma_1<1$ such that, for all $0<\gamma<\gamma_1$,
\begin{equation}
\mathbb{E}_{\theta} \|\theta^{\star} -\theta \|^2 \leq (1+h_1(\gamma)) \inf_{N \leq m_0} \bar R(\theta,N) + \frac{C_1 \epsilon^2}{n} \frac{1}{\gamma^{4\beta+1}} + \frac{C_1}{n\gamma},
\label{eq:oracle1}
\end{equation}
where
\begin{equation}
\bar R(\theta,N) = \sum_{|k|>N} |\theta_k|^2 + \frac{\epsilon^2}{n}
\sum_{|k|\leq N} |\gamma_k|^{-2} + \frac{\log^2(n)}{n} \sum_{|k|\leq N}
|\gamma_k |^{-2} |\theta_k |^2, \label{def:risque1}
\end{equation}
$h_1(\gamma) \rightarrow 0$ as $\gamma \rightarrow 0$ and $C_1$ denotes a positive constant independent of $\epsilon$ and $n$.
\end{theo}

From Theorem \ref{th:full}, our estimator $\theta^{\star}$ presents a behavior
similar to the minimizer of $\bar R(\theta,N)$. This term only differs from the
quadratic risk by a log term. This result can be explained by the choice of the
criterion (\ref{def:u1pen}). The two last terms in the right hand side of
(\ref{eq:oracle1}) are at least of order $1/n$ and may be thus considered as
negligible in most cases.

In the next section, we prove that our estimator attains the minimax of convergence on many functional spaces. In particular, the log term and the bandwidth $m_0$ have no influence on the performances of our estimator from a minimax point of view. \\

\subsubsection{Rough estimator}

In the procedure described above, we have decided to take into account the error
generated by the use of a the sequence $(\gamma_k)_{k\in\NN}$ instead of
$(\tilde \gamma_k)_{k\in \NN}$. Although their setting is slightly different
from ours, papers dealing with regularization with unknown operator consider
implicitly this error as negligible for the regularization. The goal is then to
prove that the related estimator are not affected by the noise in the operator,
i.e. this error is avoided in the oracle.

It is thus also possible to apply a similar scheme in our setting and consider
the bias enlightened in Lemma \ref{lem:decomp} as negligible. We introduce
\begin{equation}
\tilde R(\theta,N) = \sum_{|k|>N} |\theta_k|^2 + \frac{\epsilon^2}{n} \sum_{|k|\leq N} |\gamma_k|^{-2},
\label{def:risque2}
\end{equation}
that corresponds to the usual quadratic risk in an inverse problems setting.

From now on, our aim is to mimic the oracle for $\tilde R(\theta,N)$, i.e
$$ \tilde N_0 = \mathrm{arg} \min_{N\in\NN} \tilde R(\theta,N).$$
To this end, we use exactly the same scheme than for the construction of $\theta^{\star}$ starting from $\tilde R(\theta,N)$ instead of $R(\theta,N)$. Define
\begin{equation}
\tilde{U}(Y,N) = - \sum_{|k|\leq N} |\gamma_k|^{-2} \left\lbrace |\tilde c_k|^2 - \frac{\epsilon^2}{n} \right\rbrace + \frac{\epsilon^2}{n} \sum_{|k|\leq N} |\gamma_k|^{-2}.
\label{def:urisque2}
\end{equation}
Then, we introduce
\begin{equation}
\tilde N = \mathrm{arg} \min_{N\leq m_0} \tilde{U}(Y,N) \ \mathrm{and} \ \tilde \theta_k = \frac{\tilde{c}_{k}}{\gamma_{k}} \1_{\lbrace k \leq \tilde N \rbrace},
\label{def:estsimp}
\end{equation}
where $m_0$ has been introduced in (\ref{def:m0}). Hence, this estimator only
differs from the previous one by the choice of the regularization parameter
$\tilde N$. The performances of $\tilde \theta$ are detailed bellow.

\begin{theo}
\label{th:simplified} Let $\tilde \theta$ defined by (\ref{def:estsimp})  and
assume that the density $g$ satisfies Assumption \ref{assordi}. Then, there exists
$0<\gamma_2<1$ such that, for all $0<\gamma<\gamma_2$,

\begin{equation} \label{eq:oracle2}
\mathbb{E}_{\theta} \| \tilde\theta -\theta \|^2 \leq (1+h_2(\gamma)) \inf_{N \leq m_0}  R(\theta,N) +
\frac{C_2\epsilon^2}{n} \left( \frac{\|\theta\|^2\log^2(n)}{\gamma^2}
\right)^{2\beta} + \frac{C_2\epsilon^2}{n} \frac{1}{\gamma^{4\beta+1}} +
\frac{C_2}{n},
\end{equation}
 where $h_2(\gamma)\rightarrow 0$ as $\gamma \rightarrow 0$ and $C_2$ denote a positive constant independent of $\epsilon$ and $n$.
\end{theo}

We will see in Section 3 that the performances of $\theta^{\star}$ and $\tilde \theta$ are essentially the same from a minimax point of view. The existing differences may be revealed by the comparison of the oracle inequalities obtained in Theorems \ref{th:full} and \ref{th:simplified}, although this is always a difficult task. Since $\bar R(\theta,N)$ only differs from $R(\theta,N)$ by a log term, we may be interested in the residual of order $\| \theta \|^2$. For fixed $\epsilon$ and $n$, this term may have importance compared to $R(\theta,N)$, in particular for large $\| \theta \|^2$. Hence, the second estimator may be incongruous when estimating function with large norm.

More carefully, $\tilde \theta$ is a pertinent choice as soon as $\tilde
R(\theta,N)$ is close to $R(\theta,N)$. This can be strengthened by the study
of the quadratic risk defined in Lemma \ref{lem:decomp}. For instance, with a
fixed $\epsilon$, this will be the case for function with 'small' Fourier
coefficients (in particular small norms). On the other hand, as soon as
$\epsilon$ becomes 'small', the behaviour of $\tilde R(\theta,N)$ and
$R(\theta,N)$ may strongly differs. This may produce significant differences on
the performances of both $\theta^{\star}$ and $\tilde \theta$.


\section{Minimax rates of convergence for Sobolev balls} \label{sec:minimax}

We provide in this section a short discussion about the performances of our
estimator from the asymptotic minimax  point of view. For this, let $1
\leq p,q \leq \infty$  and $A>0$, and suppose that $f$ belongs to a Besov
ball $\mathcal{B}^s_{p,q}(A)$ of radius $A$ (see e.g.  \cite{DJKP95jrssb} for a precise definition of Besov spaces). \cite{BG} have derived the following asymptotic minimax lower bound for the quadratic risk over a large class of Besov balls.

\begin{theo}\label{th:minimax}
Let $1 \leq p,q \leq \infty$ and $A>0$, let $p'=p \wedge 2$ and assume that:
\begin{itemize}
\item (Regularity condition on $f$) $f \in\mathcal{B}^s_{p,q}(A)$ and $s \geq
p'$,
\item (Regularity condition on $g$) $g$ satisfies the polynomial decay condition
(\ref{eq:polydecay})  at rate $\beta$ for its Fourier coefficients,
\item (Dense case) $s \geq (2 \beta +1) (1/p - 1/2)$ and $s \geq 2\beta +1$.
\end{itemize}
Then, there exists a universal constant $M_1$ depending on $A,s,p,q$ such that
$$\inf_{\hat{f}_n} \sup_{f \in \mathcal{B}^s_{p,q}(A)} \mathbb{E} \|\hat{f}_n-f\|^2 \geq
M_1 n^{\frac{-2s }{2s + 2 \beta +1}}, \quad \text{as}\quad  n \to \infty,$$
where  $\hat{f}_{n} \in L^{2}_{p}([0,1])$ denotes any estimator of the common shape $f$, i.e a measurable function of the random processes $Y_{j}, \; j=1,\ldots,n$
\end{theo}
Therefore, Theorem \ref{th:minimax} extends the lower bound $n^{\frac{-2s }{2s + 2 \beta
+1}}$ usually obtained in a classical deconvolution model to
the more complicated model of deconvolution  with a random operator derived from
equation (\ref{eq:model-inverse}). Then, let us introduce the following smoothness class of functions which can be identified with a periodic Sobolev ball:
$$
H_{s}(A) = \left\{f \in L_{p}^{2}([0,1]) \; ; \sum_{k \in \ZZ}  (1+|k|^{2s}) |\theta_{k}|^{2} \leq A \right\},
$$
for some constant $A > 0$ and some smoothness parameter $s > 0$, where $\theta_{k} = \int_{0}^1 e^{- 2 i k \pi x} f(x) dx$. It is known (see e.g. \cite{DJKP95jrssb}) that if $s$ is not an integer then $H_{s}(A)$ can be identified with a Besov ball $\mathcal{B}^s_{2,2}(A')$.  Assuming $f \in H_{s}(A)$ with $s > 0$, then the classical choice $N^{\star} \sim n^{\frac{1}{2s+2\beta+1}}$ yields that
$$
R(\theta,N^\star)  \sim  \inf_{N \leq m_0}R(\theta,N) \sim
n^{\frac{-2s}{2s+2\beta+1}}.
$$
provided $N^{\star} \leq m_0$. It can be checked that the choice (\ref{def:m0}) implies that $m_0 \sim n^{\frac{1}{2 \beta}}$ and  thus for a sufficiently large $n$, we have that $N^\star < m_0$. Similarly the choice $N^{\star} \sim n^{\frac{1}{2s+2\beta+1}}$ yields that
$$
\bar{R}(\theta,N^\star) \sim \inf_{N \leq m_0} \bar{R}(\theta,N^\star) \sim \log^{2}(n) n
 ^{\frac{-2s}{2s+2\beta+1}},
$$
Now, remark that for the two estimators $\theta^{\star}$ and $\tilde{\theta}$,
both Theorems \ref{th:full} and \ref{th:simplified} yield that
 $\mathbb{E}_{\theta}\|\theta^{\star}-\theta\|^2 = \opO \left( \inf_{N \leq m_0}
\bar{R}(\theta,N) \right)$ and $\mathbb{E}_{\theta}\|\tilde{\theta}-\theta\|^2
= \opO \left( \inf_{N \leq m_0}R(\theta,N) \right)$ as $n \to + \infty$, since additional terms in
bounds (\ref{eq:oracle1}) and (\ref{eq:oracle2}) are of the order $\opO (\frac{1}{n^{1-\zeta}})$
for a sufficiently small positive $\zeta$. Hence, combining the above arguments one finally obtains the following result:
\begin{coro}
Suppose that the density $g$ satisfies the polynomial decay condition
(\ref{eq:polydecay})  at rate $\beta$ for its Fourier coefficients. Then, as $n \to + \infty$
$$
\sup_{f \in H_{s}(A) } \mathbb{E}_{\theta}\|\theta^{\star}-\theta\|^2 \sim \log^{2}(n) n
 ^{\frac{-2s}{2s+2\beta+1}}
$$
and
$$
\sup_{f \in H_{s}(A) } \mathbb{E}_{\theta}\|\tilde{\theta}-\theta\|^2 \sim n
 ^{\frac{-2s}{2s+2\beta+1}}.
$$
\end{coro}
From the lower bound obtained in Theorem \ref{th:minimax} we conclude that, for $s \geq 2 \beta +1$, the performances of the estimator $\tilde{\theta}$  are asymptotically optimal from the
minimax point of view, while the estimator  $\theta^{\star}$ is near-optimal up to a $ \log^{2}(n)$ factor.  This near-optimal rate of convergence of $\theta^{\star}$ is due to the use of the penalised criterion $\bar U(Y,N)$, see (\ref{def:u1pen}), with a penalty term involving a $\frac{\log^{2}(n)}{n}$ factor used to eliminate the term $ \frac{1}{n} \sum_{|k|\leq N} |\gamma_k|^{-4} | \left\lbrace |\tilde c_k|^2 -\frac{\epsilon^2}{n} \right\rbrace$ in the unbiased risk $U(Y,N)$, see (\ref{def:urisque1}). This shows that the performances of $\theta^{\star}$ and $\tilde \theta$ are essentially the same from a minimax point of view. 

\section{Numerical experiments} \label{sec:simu}

\begin{figure}[h!] 
\centering
\subfigure[]
{ \includegraphics[width=5cm]{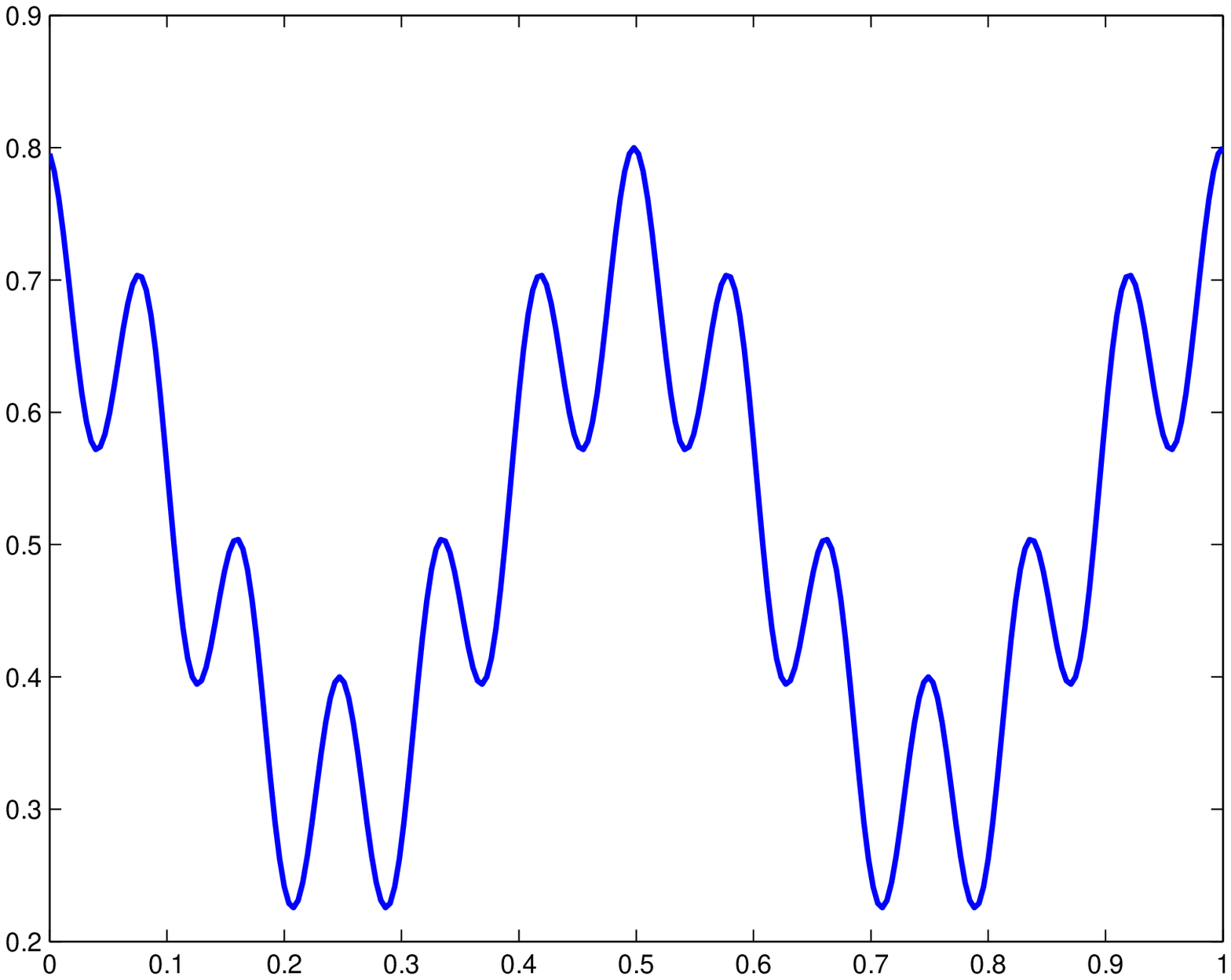} }
\subfigure[]
{ \includegraphics[width=5cm]{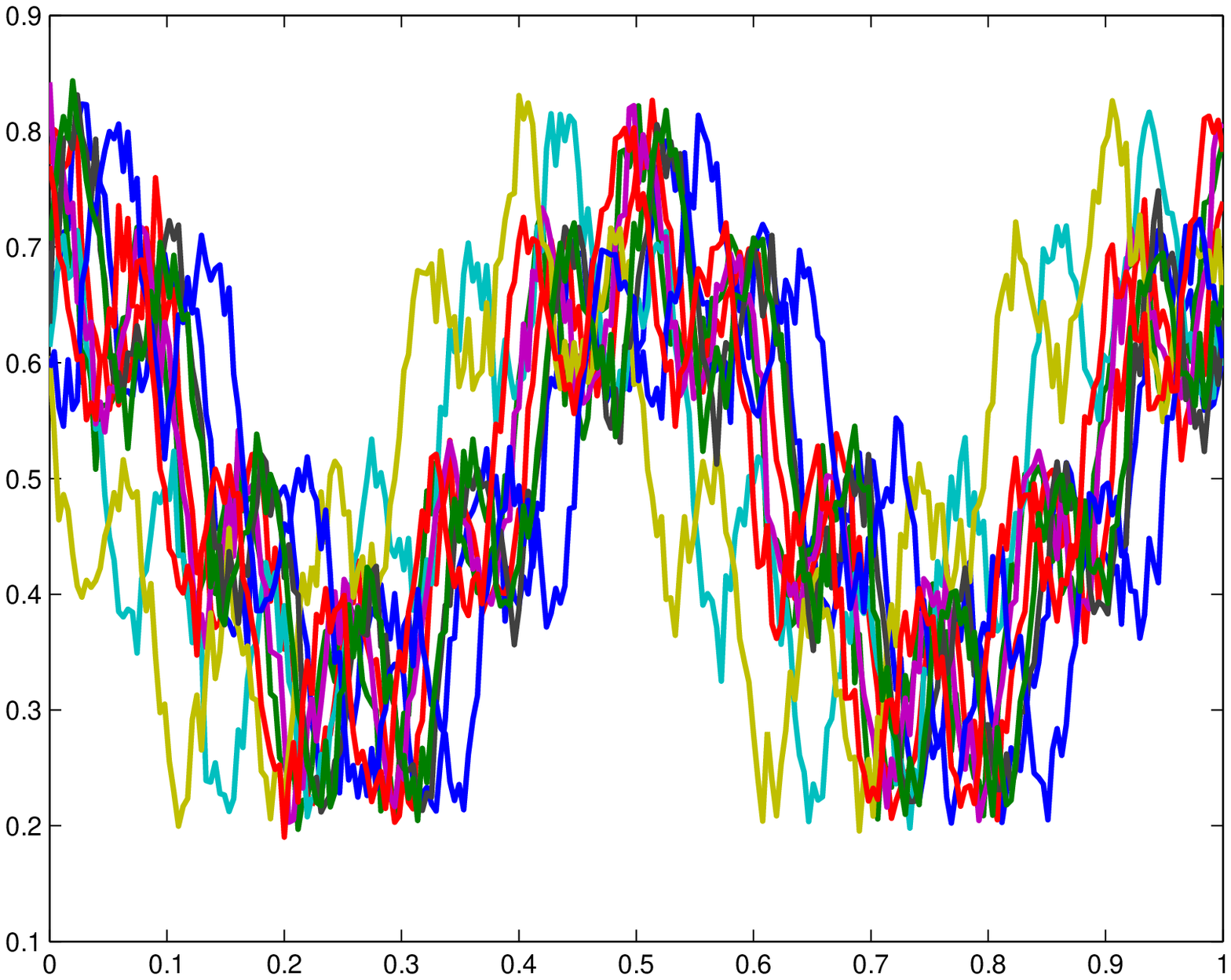} }

\subfigure[]
{ \includegraphics[width=5cm]{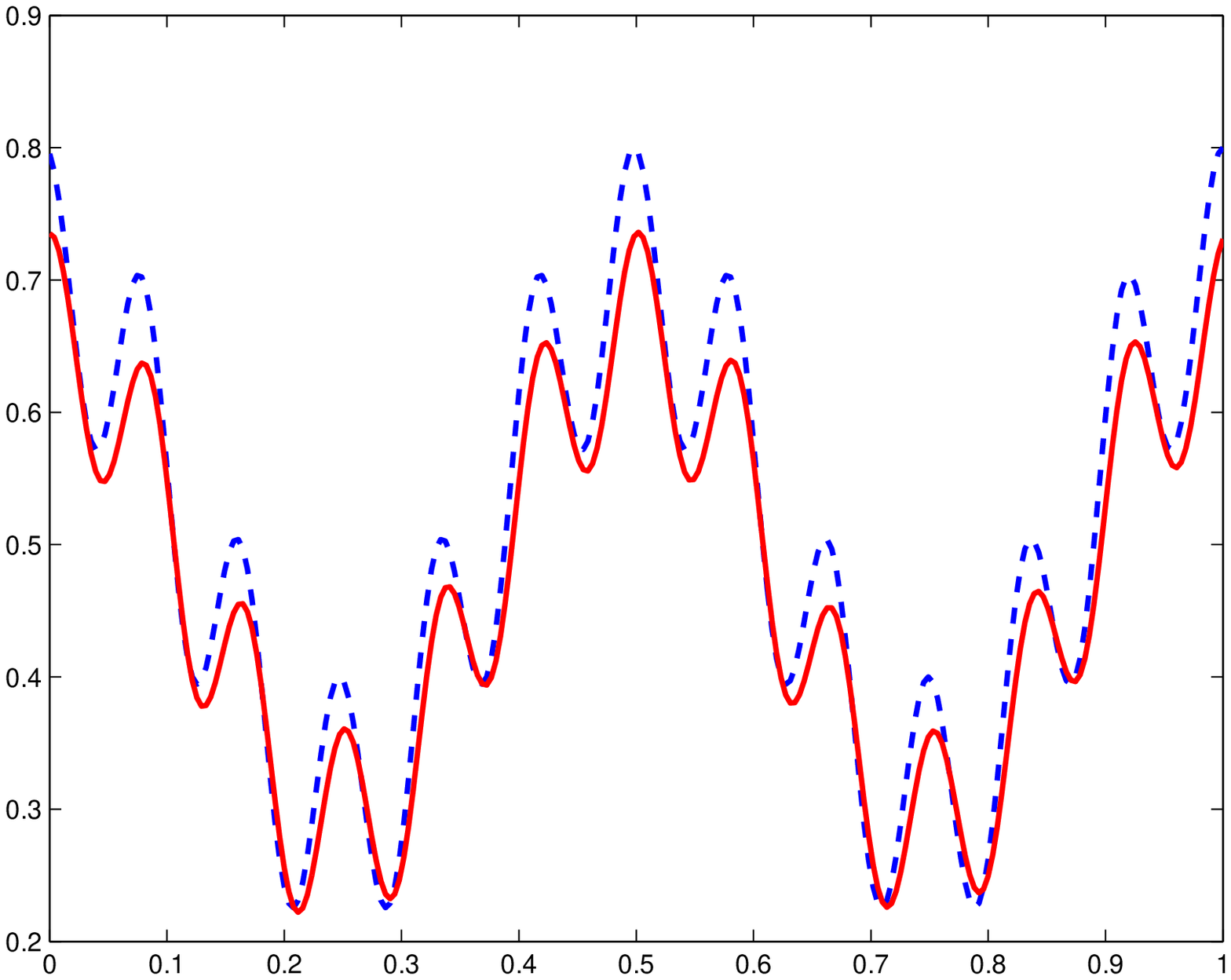} }
\subfigure[]
{ \includegraphics[width=5cm]{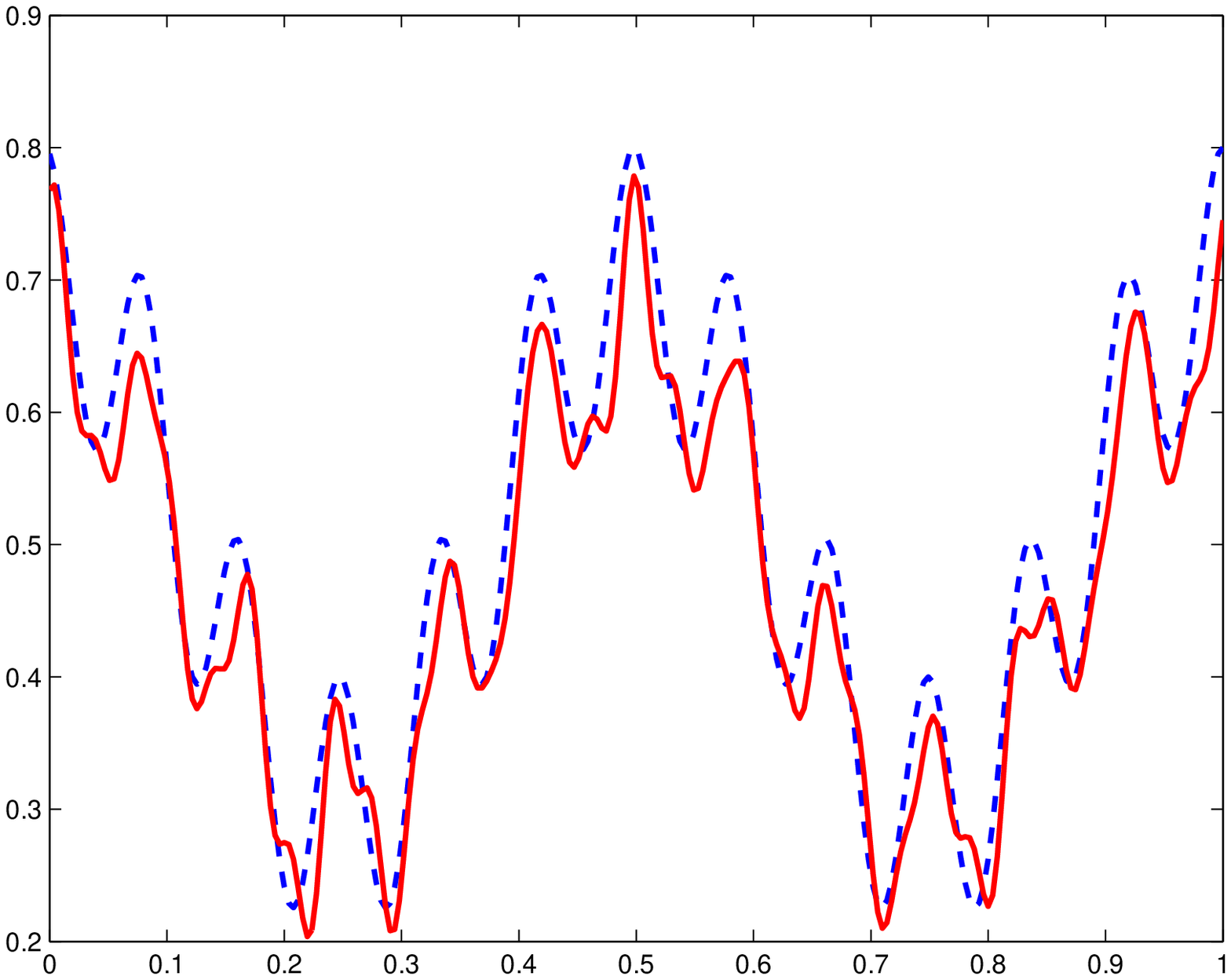} }

\caption{Wave function. (a) Mean pattern  $f$, (b) Sample of 10 curves out of $n=100$, (c) Estimation by spectral cut-off with $N^{\star} = 13$,  (d) Estimation by spectral cut-off with $\tilde{N} = 30$. The dotted curve corresponds to the true mean pattern $f$.} \label{Fig:wave}
\end{figure}

For the mean pattern $f$ to recover, we consider the smooth function shown in Figure \ref{Fig:wave}(a). Then, we simulate $n = 100$ randomly shifted curves with shifts following a Laplace distribution $g(x) = \frac{1}{\sqrt{2}\sigma} \exp \left( -\sqrt{2}\frac{|x|}{\sigma}\right)$ with $\sigma = 0.1$. Gaussian noise with a moderate variance (different to that used in the Laplace distribution) is then added to each curve. A subsample of 10 curves is shown in Figure  \ref{Fig:wave}(b). The Fourier coefficients of the density $g$ are given by $\gamma_{k} = \frac{1}{1+2\sigma^{2} \pi^{2} k^{2}}$ which corresponds to a degree of ill-posedness $\beta=2$.

The condition (\ref{def:m0}) thus leads to the choice $m_{0} = 32$. Minimisation of the criterions (\ref{def:fenetre}) and   (\ref{def:estsimp}) leads respectively to the choices $N^{\star} = 13$  and $\tilde{N} = 30$. An example of estimation by spectral cut-off using either the value of  $N^{\star}$ or $\tilde{N} $ is displayed in Figure \ref{Fig:wave}(c) and Figure \ref{Fig:wave}(d). The estimator obtained with the frequency cut-off $N^{\star} = 13$ is very satisfactory, while the choice  $\tilde{N} = 30$ seems to be too large as the resulting estimator in Figure \ref{Fig:wave}(d) is not as smooth as the estimator with  $N^{\star} = 13$.

\begin{figure}[h!]
\centering
\subfigure[]
{ \includegraphics[width=5cm]{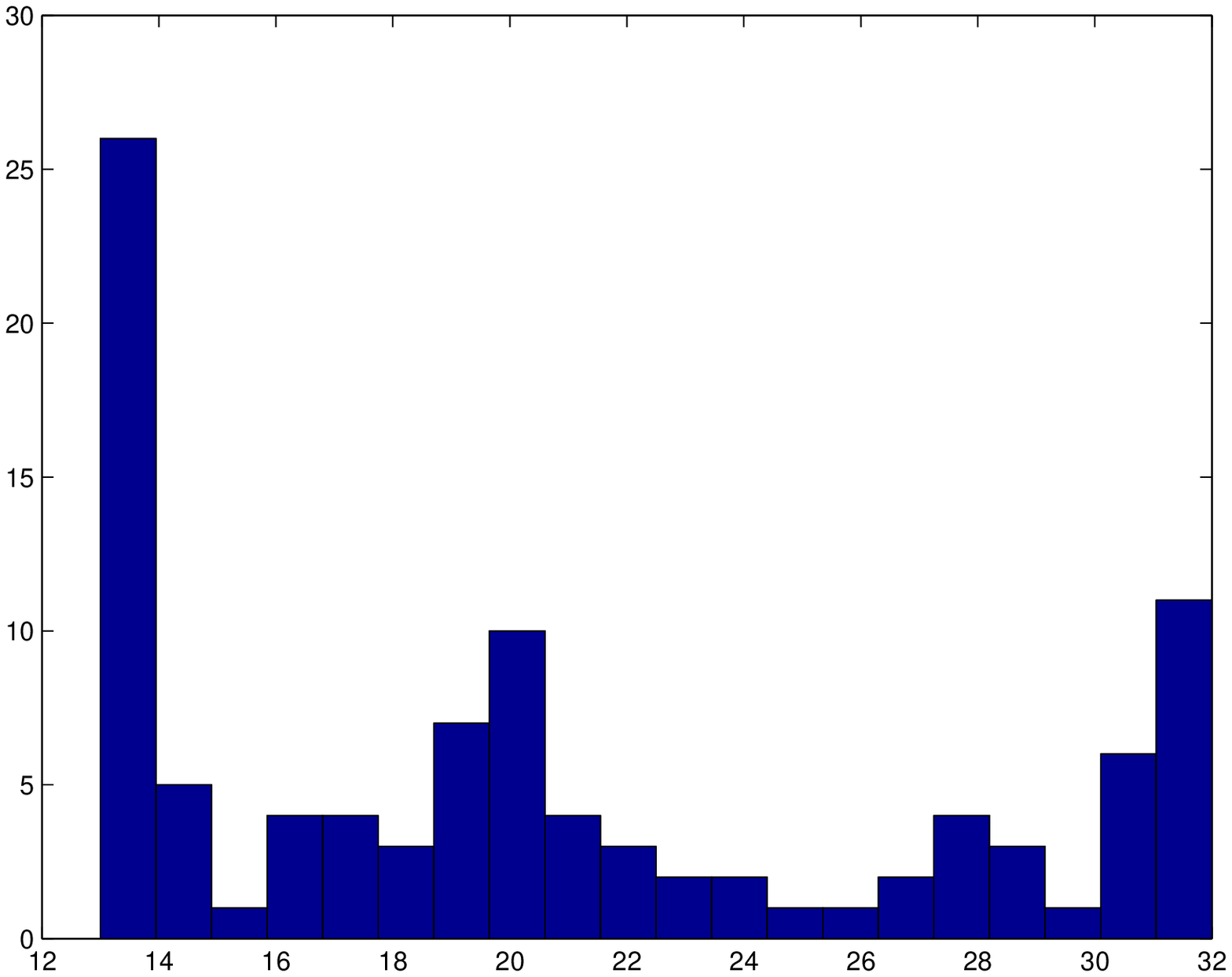} }
\subfigure[]
{ \includegraphics[width=5cm]{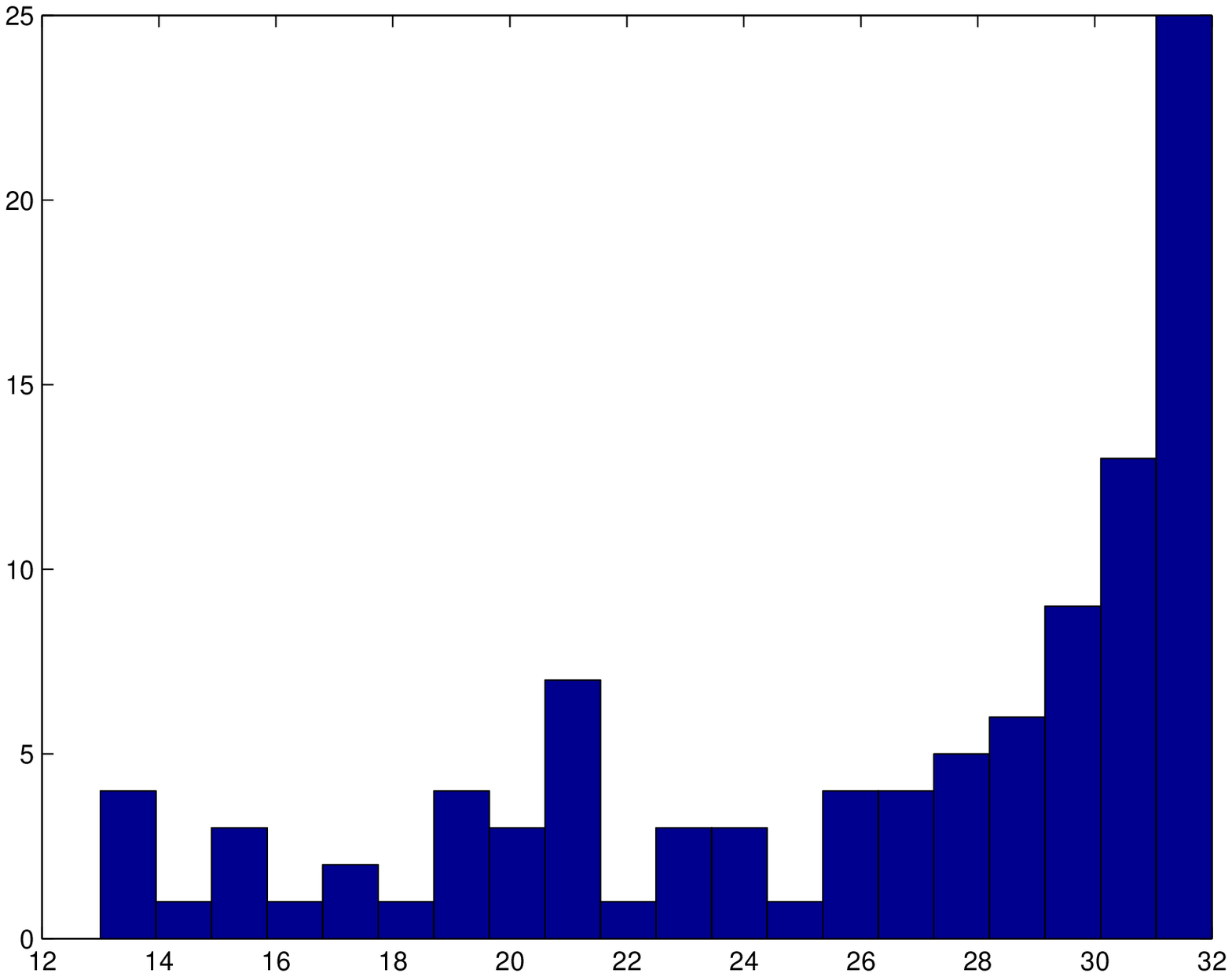} }

\caption{Selection of the frequency cut-off over $M=100$ replications of the simulations (with $m_{0} = 32$): (a) Histogram of the selected value for $N^{\star}$, (b) Histogram of the selected value for $\tilde{N} $.}  \label{Fig:simus}
\end{figure}

 This result tends to suggest that minimising $\bar U(Y,N)$ leads to a smaller choice for the frequency cut-off than the one obtained by the minimisation of the criterion $ \tilde{U}(Y,N)$. This is confirmed by the results displayed in Figure \ref{Fig:simus} which gives the histogram of the selected values for $N^{\star}$ and $\tilde{N} $ over $M=100$ independent replications of the above described simulations. Clearly the value of $N^{\star}$ is generally much smaller than $\tilde{N}$, and thus minimising  (\ref{def:estsimp}) may lead to undersmoothing which illustrates numerically our discussion in Section \ref{sec:estimf} on the differences between $\theta^{\star}$ and $\tilde{\theta}$.
 
\section{Proofs}

\noindent
\textbf{Proof of Theorem \ref{th:full}}. The proof uses the following scheme. In a first time, we compute the quadratic risk of $\theta^{\star}$ and we prove that it is close to $\bar R(\theta,N^{\star})$. The aim of the second part is to prove that $\bar U(Y,N^{\star})$ is close to $\bar R(\theta,N^{\star})$, even for a random bandwidth $N^{\star}$. Then, we use the fact that $N^{\star}$ minimizes the criterion $\bar U(Y,N^{\star})$ over the integer smaller than $m_0$ and we compute the expectation of $U(Y,N)$ for all deterministic $N$ in order to obtain an oracle inequality.

In a first time,
\begin{eqnarray*}
\E \| \theta^{\star} - \theta \|^2 & = & \E \sum_{k\in \mathbb{Z}} | \theta_k^{\star} -\theta_k |^2, \\
& = & \E \sum_{|k|\leq N^{\star}} |\gamma_k^{-1} \tilde c_k - \theta_k |^2 + \E \sum_{|k|>N^{\star}} |\theta_k|^2, \\
& = & \E \sum_{|k|\leq N^{\star}} \left| \frac{\tilde\gamma_k}{\gamma_k} \theta_k - \theta_k + \gamma_k^{-1} \frac{\epsilon}{\sqrt{n}} \xi_k \right|^2 + \E \sum_{|k|>N^{\star}} |\theta_k|^2, \\
& = & \E \sum_{|k|\leq N^{\star}} \left| \frac{\tilde\gamma_k}{\gamma_k}-1 \right|^2 |\theta_k |^2 + \frac{\epsilon^2}{n} \E \sum_{|k|\leq N^{\star}} |\xi_k|^2 |\gamma_k|^{-2}\\
& & \hspace{2cm} + \E \sum_{|k|>N^{\star}} |\theta_k|^2 + 2 \E \sum_{|k|\leq N^{\star}} \frac{\epsilon}{\sqrt{n}} Re\left( (\gamma_k^{-1}\tilde\gamma_k -1 )\theta_k \times \bar \gamma_k^{-1} \bar \xi_k \right),
\end{eqnarray*}
where for a given $z\in\mathbb{C}$, $Re(z)$ denotes the real part of $z$ and
$\bar z$ the conjuguate. The last equality can be rewritten as
\begin{eqnarray}
\E \| \theta^{\star} - \theta \|^2 & = & \E \tilde R(\theta,N^{\star}) + \E \sum_{|k|\leq N^{\star}} \left| \frac{\tilde\gamma_k}{\gamma_k}-1 \right|^2 |\theta_k |^2 + \frac{\epsilon^2}{n} \E \sum_{|k|\leq N^{\star}} |\gamma_k|^{-2} (|\xi_k|^2 -1) \nonumber \\
& & \hspace{3cm} + 2 \E \sum_{|k|\leq N^{\star}} \frac{\epsilon}{\sqrt{n}} Re\left( (\gamma_k^{-1}\tilde\gamma_k -1 )\theta_k \times \bar \gamma_k^{-1} \bar \xi_k \right), \nonumber \\
& = & \E \tilde R(\theta,N^{\star}) + A_1 + A_2 + A_3,
\label{eq:A1+A2+A3}
\end{eqnarray}
where $\tilde R(\theta,N)$ is defined in (\ref{def:risque2}). Thanks to Lemma \ref{lemme:1}, setting $K=1$,
\begin{equation}
A_1 = \E \sum_{|k|\leq N^{\star}} \left| \frac{\tilde\gamma_k}{\gamma_k}-1 \right|^2 |\theta_k |^2 \leq  \frac{\log^2(n)}{n} \E \sum_{|k|\leq N^{\star}} |\gamma_k|^{-2} |\theta_k |^2 + \frac{C}{n}.
\label{eq:A1}
\end{equation}
Now, consider a bound for $A_2$. For all $N\in\mathbb{N}$ set $\Sigma_N =
\sum_{|k|\leq N} |\gamma_k|^{-4}$. Then, for all $p\in ]1,2[$ and $1>\gamma>0$:
\begin{eqnarray*}
A_2 & = & \frac{\epsilon^2}{n} \E \sum_{|k|\leq N^{\star}} |\gamma_k|^{-2} (|\xi_k|^2 -1),\\
& = & \frac{\epsilon^2}{n} \E \left[ \sum_{|k|\leq N^{\star}} |\gamma_k|^{-2} (|\xi_k|^2 -1) - \gamma \sqrt{\Sigma_N^{\star}}^p \right] + \gamma \frac{\epsilon^2}{n} \E \sqrt{\Sigma_N^{\star}}^p, \\
& \leq & \frac{\epsilon^2}{n} \E \sup_N \left[ \sum_{|k|\leq N} |\gamma_k|^{-2} (|\xi_k|^2 -1) - \gamma \sqrt{\Sigma_N^{\star}}^p \right] + \gamma \frac{\epsilon^2}{n} \E \sqrt{\Sigma_N^{\star}}^p, \\
& \leq & \gamma \frac{\epsilon^2}{n} \E \sqrt{\Sigma_N^{\star}}^p + \frac{C}{\gamma^{1/(1-p)}} \frac{\epsilon^2}{n}.
\end{eqnarray*}
The last step can be derived from a Doob inequality: see for instance \cite{riskhull}. Thanks to the polynomial Assumption \ref{assordi} on the sequence $(\gamma_k)_k$ and setting $p=2\times (2\beta+1)/(4\beta+1)$, we obtain
\begin{equation}
\label{eq:A2}
A_2 = \frac{\epsilon^2}{n} \E \sum_{|k|\leq N^{\star}} |\gamma_k|^{-2} (|\xi_k|^2 -1) \leq \gamma \frac{\epsilon^2}{n} \E \sum_{|k|\leq N^{\star} } |\gamma_k|^{-2} + \frac{C}{\gamma^{4\beta+1}} \frac{\epsilon^2}{n}.
\end{equation}
Then, for all $1>B>0$, using the Cauchy-Schwarz and Young inequalities with the
bounds (\ref{eq:A1}) and (\ref{eq:A2})
\begin{eqnarray*}
A_3 & = & 2 \E \sum_{|k|\leq N^{\star}} \frac{\epsilon}{\sqrt{n}} Re\left( (\gamma_k^{-1}\tilde\gamma_k -1 )\theta_k \times \bar \gamma_k^{-1} \bar \xi_k \right), \nonumber\\
& \leq & B \frac{\epsilon^2}{n} \E \sum_{|k|\leq N^{\star}} |\gamma_k|^{-2} |\xi_k|^2 + B^{-1} \E \sum_{|k|\leq N^{\star}} |\theta_k|^2 \left| \frac{\tilde\gamma_k}{\gamma_k} -1 \right|^2, \nonumber\\
\end{eqnarray*}
Thus, for any $K>0$,
\begin{equation}A_3  \leq  (B+ B\gamma) \frac{\epsilon^2}{n} \E \sum_{|k|\leq N^{\star}}
|\gamma_k|^{-2} + B^{-1}K \frac{\log^2(n)}{n} \E \sum_{|k|\leq N^{\star}}
|\gamma_k|^{-2} |\theta_k |^2 + \frac{C\epsilon^2}{n\gamma^{4\beta+1}}+
\frac{C}{nK}. \label{eq:A3}
\end{equation}
With $B=\sqrt{K}=\sqrt{\gamma}$, we obtain from (\ref{eq:A1+A2+A3})-(\ref{eq:A3})
\begin{equation}
\E \|\theta^{\star} -\theta \|^2 \leq (1+\gamma+2\sqrt{\gamma}) \E \bar R(\theta,N^{\star}) + \frac{C\epsilon^2}{n\gamma^{4\beta+1}}+ \frac{C}{n},
\label{eq:etape1}
\end{equation}
where $\bar R(\theta,N)$ is defined in (\ref{def:risque1}). This concludes the first step of our proof. Now, we write $\bar U(Y,N^{\star})$ in terms of $\bar R(\theta,N^{\star})$. In the following, we define $x_n =(1-n^{-1})$. We have
\begin{eqnarray*}
\bar U(Y,N^{\star}) & = & -x_n \sum_{|k|\leq N^{\star}} |\gamma_k|^{-2} \left\lbrace |\tilde c_k |^2 -\frac{\epsilon^2}{n} \right\rbrace + \frac{\epsilon^2}{n} \sum_{|k|\leq N^{\star}} |\gamma_k|^{-2}+ \frac{\log^2(n)}{n} \sum_{|k|\leq N} |\gamma_k |^{-4}  \left\lbrace |\tilde c_k| -\frac{\epsilon^2}{n} \right\rbrace,\\
& = & \bar R(\theta,N^{\star}) -\left(1-\frac{1}{n} \right) \sum_{|k|\leq N^{\star}} |\gamma_k|^{-2} \left\lbrace |\tilde c_k |^2 -\frac{\epsilon^2}{n} \right\rbrace - \sum_{|k|\geq N^{\star}} |\theta_k|^2\\
& & \hspace{4cm} + \frac{\log^2(n)}{n} \sum_{|k|\leq N^{\star}} \left[ |\gamma_k |^{-4}  \left\lbrace |\tilde c_k| -\frac{\epsilon^2}{n} \right\rbrace - |\gamma_k |^{-2} |\theta_k |^2 \right]
\end{eqnarray*}
This equality can be rewritten as
\begin{eqnarray}
\bar R(\theta,N^{\star}) & = & \bar U(Y,N^{\star}) + \|\theta \|^2 + x_n\sum_{|k|\leq N^{\star}} \left\lbrace |\gamma_k|^{-2} |\tilde c_k |^2 -\frac{\epsilon^2}{n} |\gamma_k |^{-2} \right\rbrace - \sum_{|k|\leq N^{\star}} |\theta_k|^2 \nonumber \\
& & \hspace{3cm} + \frac{\log^2(n)}{n} \sum_{|k|\leq N^{\star}} \left[ |\gamma_k |^{-2} |\theta_k|^2 - |\gamma_k|^{-4} \left\lbrace |\tilde c_k|^2 - \frac{\epsilon^2}{n}  \right\rbrace \right].
\label{eq:inter1}
\end{eqnarray}
For all $k\in\mathbb{N}$
$$ | \tilde c_k |^2 = | \theta_k \tilde\gamma_k |^2 + \frac{\epsilon^2}{n} |\xi_k|^2 + 2 \epsilon n^{-1/2} Re(\theta_k \tilde\gamma_k \bar\xi_k),$$
and
$$|\gamma_k |^{-2} |\tilde c_k |^2 = |\theta_k|^2 \left| \frac{\tilde \gamma_k}{\gamma_k} \right|^2 + \frac{\epsilon^2}{n} |\gamma_k |^{-2} |\xi_k |^2 + 2 \frac{\epsilon}{\sqrt{n}} |\gamma_k |^{-2} Re(\theta_k \tilde\gamma_k \bar\xi_k).$$
Since $x_n<1$
\begin{eqnarray}
\lefteqn{x_n \E \sum_{|k|\leq N^{\star}} \left\lbrace |\gamma_k|^{-2} |\tilde c_k |^2 -\frac{\epsilon^2}{n} |\gamma_k |^{-2}\right\rbrace - \sum_{|k|\leq N^{\star}} |\theta_k|^2 } \nonumber \\
& \leq &  \E \sum_{|k|\leq N^{\star}} |\theta_k|^{2} \left( \left| \frac{\tilde \gamma_k}{\gamma_k} \right|^2 -1 \right)
+ \frac{\epsilon^2}{n} x_n \sum_{|k|\leq N^{\star}} |\gamma_k|^{-2} (|\xi_k|^2-1) + 2\frac{\epsilon}{\sqrt{n}} x_n \sum_{|k|\leq N^{\star}} |\gamma_k|^{-2}Re(\theta_k \tilde\gamma_k \bar\xi_k), \nonumber \\
& = &  E_1 + E_2 + E_3.
\label{eq:E1+E2+E3}
\end{eqnarray}

First consider the bound of $E_1$. Thanks
to Lemma \ref{lemme:2} and some simple algebra
\begin{eqnarray*}
E_1 & = & \E \sum_{|k|\leq N^{\star}} |\theta_k|^{2} \left( \left| \frac{\tilde \gamma_k}{\gamma_k} \right|^2 -1 \right) ,\\
& \leq &  2\gamma \frac{\log^2(n)}{n} \E \sum_{|k|\leq N^{\star}} |\theta_k|^2 |\gamma_k|^{-2} + \gamma \E \sum_{|k|> N^{\star}} |\theta_k|^2 \\
& & \hspace{1cm} + \gamma  \sum_{|k|> N_0} |\theta_k|^2 + \gamma^{-1} \sum_{|k|\leq N_0} |\theta_k|^2 |\gamma_k|^{-2}(1-|\gamma_k|^2) + \frac{C}{n \gamma^2},  \\
& \leq & 2\gamma \E \bar R(\theta,N^{\star}) + \left( \gamma + \frac{\gamma^{-1}}{\log^2(n)} \right) \bar R(\theta,N_0) + \frac{C}{n\gamma^2},
\end{eqnarray*}
where
$$N_0 = \mathrm{arg} \min_{N\leq m_0} \bar R(\theta,N).$$
The terms $E_2$ and $E_3$ are bounded using respectively (\ref{eq:A2}) and Lemma \ref{lemme:3}. We get
\begin{eqnarray}
\lefteqn{\E \sum_{|k|\leq N^{\star}} \left\lbrace |\gamma_k|^{-2} |\tilde c_k |^2 -\frac{\epsilon^2}{n} |\gamma_k |^{-2}- |\theta_k|^2  \right\rbrace} \nonumber \\
& \leq & D\gamma \E \bar R(\theta,N^{\star}) + D\gamma \bar R(\theta,N_0)  + \frac{\epsilon^2}{n}\frac{C}{\gamma^{4\beta+1}} + \frac{C}{n\gamma^2}.
\label{eq:numm1}
\end{eqnarray}
We are now interested in the second residual term of (\ref{eq:inter1}). Thanks to the definition of $\tilde c_k$:
\begin{eqnarray}
\lefteqn{\frac{\log^2 n}{n} \E \sum_{|k|\leq N^{\star}} |\gamma_k|^{-2} \left\lbrace -|\gamma_k|^{-2} |\tilde c_k |^2 +\frac{\epsilon^2}{n} |\gamma_k |^{-2} + |\theta_k|^2  \right\rbrace} \nonumber \\
& = &  \E \sum_{|k|\leq N^{\star}} |\gamma_k|^{-2}|\theta_k|^{2} \left( 1 - \left| \frac{\tilde \gamma_k}{\gamma_k} \right|^2  \right)
+ \frac{\epsilon^2}{n} \sum_{|k|\leq N^{\star}} |\gamma_k|^{-4} (1 - |\xi_k|^2) - 2\frac{\epsilon}{\sqrt{n}} \sum_{|k|\leq N^{\star}} |\gamma_k|^{-4}Re(\theta_k \tilde\gamma_k \bar\xi_k), \nonumber \\
& \leq & D\gamma \E \bar R(\theta,N^{\star}) + D\gamma \bar R(\theta,N_0)  + \frac{\epsilon^2}{n}\frac{C}{\gamma^{4\beta+1}} + \frac{C}{n\gamma^2},
\label{eq:numm2}
\end{eqnarray}
for some $D>0$ independent of $\epsilon$ and $n$.Indeed, we can use essentialy the same algebra as for the bound of the terms $E_1, E_2$ and $E_3$ and the inequality
$$ |\gamma_k |^{-2} \leq \frac{n}{\log^2n}, \ \forall k\leq m_0.$$
Hence, using (\ref{eq:numm1}) and (\ref{eq:numm2})
\begin{equation}
(1-D\gamma)\E\bar R(\theta,N^{\star}) \leq \E U(Y,N^{\star}) + \|\theta\|^2 + D\gamma \tilde R(\theta,N_0) + \frac{C}{n\gamma^2} + \frac{C\epsilon^2}{n} \frac{1}{\gamma^{4\beta+1}}.
\label{eq:etape2}
\end{equation}
>From the definition of $N^{\star}$, we immediatly get
$$(1-D\gamma)\E\bar R(\theta,N^{\star}) \leq \E U(Y,N_0) + \|\theta\|^2 + D\gamma \tilde R(\theta,N_0) + \frac{C}{n\gamma^2} + \frac{C\epsilon^2}{n} \frac{1}{\gamma^{4\beta+1}},$$
where $N_0$ denotes the oracle bandwidth. Since
$$\E U(Y,N_0) = \tilde R(\theta,N_0) - \|\theta\|^2,$$
we obtain
\begin{equation}
(1-D\gamma) \E \tilde R(\theta,N) \leq (1+D\gamma) \tilde R(\theta,N_0) + \frac{C}{n\gamma^2} + \frac{C\epsilon^2}{n} \frac{1}{\gamma^{4\beta+1}}.
\label{eq:etapefin}
\end{equation}
Using (\ref{eq:etape1}) and (\ref{eq:etapefin}), we get:
\begin{eqnarray*}
\E\|\theta^{\star} -\theta \|^2
& \leq & (1+D\sqrt{\gamma}) \E \tilde R(\theta,N^{\star}) + \frac{C\epsilon^2}{n} \frac{1}{\gamma^{4\beta+1}} + \frac{C}{n\gamma},\\
& \leq & \left( \frac{1+D\sqrt{\gamma}}{1-D\gamma} \right) \tilde R(\theta,N_0) + \frac{C\epsilon^2}{n} \frac{1}{\gamma^{4\beta+1}} + \frac{C}{n\gamma}.
\end{eqnarray*}
This concludes the proof of Theorem \ref{th:full}.
\begin{flushright}
$\Box$
\end{flushright}

\noindent
\textbf{Proof of Theorem \ref{th:simplified}}. The proof follows the same main lines as for Theorem \ref{th:full}. Inequality (\ref{eq:A1+A2+A3}) provides:
\begin{eqnarray*}
\E \| \theta^{\star} - \theta \|^2 & = & \E \tilde R(\theta,N^{\star}) + \E \sum_{|k|\leq N^{\star}} \left| \frac{\tilde\gamma_k}{\gamma_k}-1 \right|^2 |\theta_k |^2 + \frac{\epsilon^2}{n} \E \sum_{|k|\leq N^{\star}} |\gamma_k|^{-2} (|\xi_k|^2 -1) \\
& & \hspace{3cm} + 2 \E \sum_{|k|\leq N^{\star}} \frac{\epsilon}{\sqrt{n}} Re\left( (\gamma_k^{-1}\tilde\gamma_k -1 )\theta_k \times \bar \gamma_k^{-1} \bar \xi_k \right), \\
& = & \E \tilde R(\theta,N^{\star}) + A_1 + A_2 + A_3.
\end{eqnarray*}
Thanks to Lemma \ref{lemme:1} and an inequality of \cite{cavgopitsy}, we obtain
for all $0<\gamma<1$:
\begin{eqnarray}
A_1 & \leq & \log^2(n) \frac{\epsilon^2}{n} \E \sup_{|k|\leq N^{\star}} |\gamma_k|^{-2} |\theta_k |^2 + \frac{C}{n} \nonumber,\\
& \leq & \gamma \frac{\epsilon^2}{n} \E \sum_{|k|\leq N^{\star}} |\gamma_k|^{-2} + \frac{C\epsilon^2}{n} \left( \frac{\|\theta\|^2 \log^2 (n)}{\gamma} \right)^{2\beta} + \frac{C}{n}.
\label{eq:A1'}
\end{eqnarray}
Then, for all $B>0$, using the Cauchy-Schwarz and Young inequalities with the
bounds (\ref{eq:A1}),(\ref{eq:A2}):
\begin{eqnarray}
A_3 & = & 2 \E \sum_{|k|\leq N^{\star}} \frac{\epsilon}{\sqrt{n}} Re\left( (\gamma_k^{-1}\tilde\gamma_k -1 )\theta_k \times \bar \gamma_k^{-1} \bar \xi_k \right), \nonumber\\
& \leq & (B+ B\gamma + B^{-1}\gamma) \frac{\epsilon^2}{n} \E \sum_{|k|\leq N^{\star}} |\gamma_k|^{-2} + \frac{C\epsilon^2}{n} \left( \frac{\|\theta\|^2 \log^2 (n)}{\gamma} \right)^{2\beta} + \frac{C\epsilon^2}{n\gamma^{4\beta+1}}+ \frac{C}{n}.
\label{eq:A3'}
\end{eqnarray}
With the choice $B=\sqrt{\gamma}$, we obtain from (\ref{eq:A1+A2+A3})-(\ref{eq:A3}):
\begin{equation}
\E \|\theta^{\star} -\theta \|^2 \leq (1+3\gamma+2\sqrt{\gamma}) \E \tilde R(\theta,N^{\star}) + \frac{C\epsilon^2}{n} \left( \frac{\|\theta\|^2 \log^2 (n)}{\gamma} \right)^{2\beta} + \frac{C\epsilon^2}{n\gamma^{4\beta+1}}+ \frac{C}{n}.
\label{eq:etape1'}
\end{equation}
Then,
\begin{eqnarray*}
U(Y,N^{\star}) & = & -\sum_{|k|\leq N^{\star}} \left\lbrace |\gamma_k|^{-2} |\tilde c_k |^2 -\frac{\epsilon^2}{n} |\gamma_k |^{-2} \right\rbrace + \frac{\epsilon^2}{n} \sum_{|k|\leq N^{\star}} |\gamma_k|^{-2},\\
& = & -\sum_{|k|\leq N^{\star}} \left\lbrace |\gamma_k|^{-2} |\tilde c_k |^2 -\frac{\epsilon^2}{n} |\gamma_k |^{-2} \right\rbrace - \sum_{|k|\geq N^{\star}} |\theta_k|^2 + \sum_{|k|\geq N^{\star}} |\theta_k|^2 + \frac{\epsilon^2}{n} \sum_{|k|\leq N^{\star}} |\gamma_k|^{-2},\\
& = & \tilde R(\theta,N^{\star}) -\sum_{|k|\leq N^{\star}} \left\lbrace |\gamma_k|^{-2} |\tilde c_k |^2 -\frac{\epsilon^2}{n} |\gamma_k |^{-2} \right\rbrace - \sum_{|k|\geq N^{\star}} |\theta_k|^2.
\end{eqnarray*}
This equality can be rewritten as
$$\tilde R(\theta,N^{\star}) = U(Y,N^{\star}) + \|\theta \|^2 + \sum_{|k|\leq N^{\star}} \left\lbrace |\gamma_k|^{-2} |\tilde c_k |^2 -\frac{\epsilon^2}{n} |\gamma_k |^{-2}- |\theta_k|^2  \right\rbrace.$$
Hence,
\begin{eqnarray}
\E \tilde R(\theta,N) & = & \E U(Y,N^{\star}) + \|\theta \|^2 + \E \sum_{|k|\leq N^{\star}} |\theta_k|^{2} \left( \left| \frac{\tilde \gamma_k}{\gamma_k} \right|^2 -1 \right) \nonumber \\
& & \hspace{2cm} + \frac{\epsilon^2}{n} \sum_{|k|\leq N^{\star}} |\gamma_k|^{-2} (|\xi_k|^2-1) + 2\frac{\epsilon}{\sqrt{n}} \sum_{|k|\leq N^{\star}} |\gamma_k|^{-2}Re(\theta_k \tilde\gamma_k \bar\xi_k), \nonumber \\
& = & \E U(Y,N^{\star}) + \|\theta \|^2 + E_1 + E_2 + E_3.
\end{eqnarray}
Using previous results:
\begin{eqnarray}
E_1 & = & \E \sum_{|k|\leq N^{\star}} |\theta_k|^{2} \left( \left| \frac{\tilde \gamma_k}{\gamma_k} \right|^2 -1 \right) , \nonumber\\
& \leq & 2\gamma \E \tilde R(\theta,N^{\star}) + \gamma \tilde R(\theta,N_0) + \frac{C}{n} + \frac{C\epsilon^2}{n} \left( \frac{\|\theta\|^2 \log^2(n)}{\gamma^2} \right)^{2\beta}.
\label{eq:E1'}
\end{eqnarray}
The terms $E_2$ and $E_3$ are bounded using respectively (\ref{eq:A2}) and Lemma \ref{lemme:3}. We get:
\begin{eqnarray*}
\E \tilde R(\theta,N^{\star})
& \leq & \E U(Y,N^{\star}) + \|\theta\|^2 + D\gamma \E \tilde R(\theta,N^{\star}) + D\gamma \tilde R(\theta,N_0) \\
&  & + \frac{C}{n\gamma^{2\beta}} + \frac{C\epsilon^2}{n} \left( \frac{\|\theta\|^2 \log^2(n)}{\gamma^2} \right)^{2\beta}.
\end{eqnarray*}
Hence,
\begin{equation}
(1-D\gamma)\E\tilde R(\theta,N^{\star}) \leq \E U(Y,N^{\star}) + \|\theta\|^2 + D\gamma \tilde R(\theta,N_0) + \frac{C}{n\gamma^{2\beta}} + \frac{C\epsilon^2}{n} \left( \frac{\|\theta\|^2 \log^2(n)}{\gamma^2} \right)^{2\beta}.
\end{equation}
>From the definition of $N^{\star}$, we immediatly get:
\begin{equation}
(1-D\gamma)\E\tilde R(\theta,N^{\star}) \leq \E U(Y,N_0) + \|\theta\|^2 + D\gamma \tilde R(\theta,N_0) + \frac{C}{n\gamma^{2\beta}} + \frac{C\epsilon^2}{n} \left( \frac{\|\theta\|^2 \log^2(n)}{\gamma^2} \right)^{2\beta}.
\label{eq:etape3}
\end{equation}
In order to conclude the proof, we prove that $\E U(Y,N_0)$ is close to $R(\theta,N_0)$. First remark that:
\begin{eqnarray*}
\E U(Y,N_0) & = & \E \left[ -\sum_{k=1}^{N_0} |\gamma_k|^{-2} \left\lbrace |\tilde c_k |^2 - \frac{\epsilon^2}{n} \right\rbrace  + \frac{\epsilon^2}{n} \sum_{k=1}^{N_0} |\gamma_k|^{-2} \right], \\
& = & \E \left[ -\sum_{k=1}^{N_0} \left\lbrace |\gamma_k|^{-2} |\tilde c_k |^2 - |\gamma_k|^{-2} \frac{\epsilon^2}{n} - |\theta_k|^2 \right\rbrace \right] - \sum_{k=1}^{N_0} |\theta_k|^2 + \frac{\epsilon^2}{n} \sum_{k=1}^{N_0} |\gamma_k|^{-2}.
\end{eqnarray*}
Since for all $k\in\mathbb{N}$:
$$\E |\tilde c_k |^2 = |\theta_k|^2 \E |\tilde \gamma_k|^2 + \frac{\epsilon^2}{n} = |\theta_k|^2 \left( |\gamma_k|^2 + \frac{1}{n} \right) + \frac{\epsilon^2}{n},$$
we obtain,
$$ \E U(Y,N_0) = -\sum_{k=1}^{N_0} |\theta_k|^2 \frac{|\gamma_k|^{-2}}{n} + \sum_{|k|>N_0} |\theta_k|^2 + \frac{\epsilon^2}{n} \sum_{k=1}^{N_0} |\gamma_k|^{-2} - \| \theta \|^2.$$
Therefore,
$$\E U(Y,N_0) = -\sum_{k=1}^{N_0} |\theta_k|^2 \frac{|\gamma_k|^{-2}}{n} + \tilde R(\theta,N_0) -\|\theta\|^2 \leq \tilde R(\theta,N_0) - \|\theta\|^2,$$
and
\begin{equation}
(1-D\gamma) \E \tilde R(\theta,N) \leq (1+D\gamma) \tilde R(\theta,N_0) + \frac{C}{n\gamma^{2\beta}} + \frac{C\epsilon^2}{n} \left( \frac{\|\theta\|^2\log^2(n)}{\gamma^2} \right)^{2\beta}.
\label{eq:etape4}
\end{equation}
Using (\ref{eq:etape1}) and (\ref{eq:etape4}), we get:
\begin{eqnarray*}
\E\|\theta^{\star} -\theta \|^2
& \leq & (1+D\sqrt{\gamma}) \E \tilde R(\theta,N^{\star}) + \frac{C\epsilon^2}{n} \left( \frac{\|\theta\|^2\log^2(n)}{\gamma^2} \right)^{2\beta} + \frac{C\epsilon^2}{n} \frac{1}{\gamma^{4\beta+1}} + \frac{1}{n},\\
& \leq & \left( \frac{1+D\sqrt{\gamma}}{1-D\gamma} \right) \tilde R(\theta,N_0) + \frac{C\epsilon^2}{n} \left( \frac{\|\theta\|^2\log^2(n)}{\gamma^2} \right)^{2\beta} + \frac{C\epsilon^2}{n} \frac{1}{\gamma^{4\beta+1}} + \frac{1}{n}.
\end{eqnarray*}
Since $\tilde R(\theta,N) \leq R(\theta,N)$, we eventually get:
$$ \E\|\theta^{\star} -\theta \|^2 \leq \left( \frac{1+D\sqrt{\gamma}}{1-D\gamma} \right) \inf_N R(\theta,N) + \frac{C\epsilon^2}{n} \left( \frac{\|\theta\|^2\log^2(n)}{\gamma^2} \right)^{2\beta} + \frac{C\epsilon^2}{n} \frac{1}{\gamma^{4\beta+1}} + \frac{1}{n}.$$
This concludes the proof of Theorem \ref{th:simplified}.
\begin{flushright}
$\Box$
\end{flushright}

\section*{Appendix}

\begin{lemma}
\label{lemme:1}
For all $K>0$, we have
$$ \E \sum_{|k|\leq N^{\star}} \left| \frac{\tilde\gamma_k}{\gamma_k}-1 \right|^2 |\theta_k |^2 \leq K \frac{\log^2(n)}{n} \E \sum_{|k|\leq N^{\star}} |\gamma_k|^{-2} |\theta_k |^2 + \frac{C}{nK},$$
where $C$ denote a positive constant independent of $\epsilon$ and $n$.
\end{lemma}
PROOF. Let $Q>0$ a deterministic term which will be chosen later.
\begin{eqnarray*}
\lefteqn{ \E \sum_{|k|\leq N^{\star}} \left| \frac{\tilde\gamma_k}{\gamma_k}-1 \right|^2 |\theta_k |^2 = \E \sum_{|k|\leq N^{\star}}  |\theta_k |^2 |\gamma_k|^{-2} |\tilde\gamma_k - \gamma_k |^2, }\\
& \leq & Q  \E \sum_{|k|\leq N^{\star}}  |\theta_k |^2 |\gamma_k|^{-2} + \E \sum_{|k|\leq N^{\star}} |\theta_k |^2 |\gamma_k|^{-2} \left\lbrace |\tilde\gamma_k - \gamma_k |^2 - Q \right\rbrace \1_{\lbrace |\tilde\gamma_k - \gamma_k|^2 >Q \rbrace}.
\end{eqnarray*}
Thanks to (\ref{def:fenetre}) and (\ref{def:m0})
\begin{eqnarray*}
\lefteqn{\E \sum_{|k|\leq N^{\star}}  |\theta_k |^2 |\gamma_k|^{-2} \left\lbrace |\tilde\gamma_k - \gamma_k |^2 - Q \right\rbrace \1_{\lbrace |\tilde\gamma_k - \gamma_k|^2 >Q \rbrace}}\\
& \leq & C \frac{n}{\log^2(n)} \sum_{|k|\leq m_0} |\theta_k|^2 \E \left\lbrace |\tilde\gamma_k - \gamma_k |^2 - Q \right\rbrace \1_{\lbrace |\tilde\gamma_k - \gamma_k|^2 >Q \rbrace}.
\end{eqnarray*}
For all $|k|\leq m_0$, using an integration by part
$$\E \left[ |\tilde\gamma_k - \gamma_k |^2 - Q \right] \1_{\lbrace |\tilde\gamma_k - \gamma_k|^2 >Q \rbrace} = \int_Q^{+\infty} P(|\tilde\gamma_k -\gamma_k |^2 \geq x ) dx.$$
Let $x\geq Q$. A Bernstein type inequality provides
\begin{eqnarray*}
P(|\tilde\gamma_k -\gamma_k |^2 \geq x ) & = & P\left( \left| \frac{1}{n} \sum_{l=1}^n \left\lbrace e^{-2i\pi k \tau_l} - \mathbb{E}[ e^{-2i\pi k \tau_l}] \right\rbrace \right| \geq \sqrt{x} \right), \\
& \leq & 2 \exp \left\lbrace - \frac{(n\sqrt{x})^2}{2 \sum_{l=1}^n \mathrm{Var}(e^{-2i\pi k \tau_l}) + n\sqrt{x}/3 } \right\rbrace,\\
& \leq & 2 \exp \left\lbrace - \frac{(n\sqrt{x})^2}{2n+ n\sqrt{x}/3} \right\rbrace.
\end{eqnarray*}
Hence, for all $|k|\leq m_0$,
\begin{eqnarray*}
\E \left[ |\tilde\gamma_k - \gamma_k |^2 - Q \right] \1_{\lbrace |\tilde\gamma_k - \gamma_k|^2 >Q \rbrace}
& \leq & \int_Q^{+\infty} \exp\left\lbrace - \frac{nx}{2+ \sqrt{x}/3} \right\rbrace dx,\\
& \leq &  \int_Q^{36} \exp\left\lbrace - \frac{nx}{4} \right\rbrace dx + \int_{36}^{+\infty} \exp\left\lbrace -Cn\sqrt{x} \right\rbrace dx \leq \frac{C}{n} e^{-Qn/4},
\end{eqnarray*}
where $C$ denotes a positive constant independent of $Q$. Let $K>0$. Choosing for instance $Q=n^{-1} K\log^2(n)$, we obtain
\begin{eqnarray*}
\E \sum_{|k|\leq N^{\star}} \left| \frac{\tilde\gamma_k}{\gamma_k}-1 \right|^2 |\theta_k |^2
& \leq & K \frac{\log^2(n)}{n} \E \sum_{|k|\leq N^{\star}} |\gamma_k|^{-2} |\theta_k |^2 + \frac{Cn m_0}{\log^2(n)} e^{-K\log^2 (n)/4},\\
& \leq & K \frac{\log^2(n)}{n} \E \sum_{|k|\leq N^{\star}} |\gamma_k|^{-2} |\theta_k |^2 + \frac{C}{nK},\\
\end{eqnarray*}
where $C$ denotes a positive constant independent of $\epsilon$ and $n$. This concludes the proof of Lemma \ref{lemme:1}.
\begin{flushright}
$\Box$
\end{flushright}

\begin{lemma}
\label{lemme:2}
Let $N^{\star}$ defined in (\ref{def:fenetre}). For all deterministic bandwidth $N$ and $0<\gamma<1$, we have
\begin{eqnarray*}
\E \sum_{|k|\leq N^{\star}} |\theta_k|^{2} \left( \left| \frac{\tilde \gamma_k}{\gamma_k} \right|^2 -1 \right)
& \leq &  2\gamma \frac{\log^2(n)}{n} \E \sum_{|k|\leq N^{\star}} |\theta_k|^2 |\gamma_k|^{-2} + \gamma \E \sum_{|k|> N^{\star}} |\theta_k|^2 \\
& & \hspace{1cm} + \gamma  \sum_{|k|> N} |\theta_k|^2 + \frac{\gamma^{-1}}{n} \sum_{|k|\leq N} |\theta_k|^2 |\gamma_k|^{-2}(1-|\gamma_k|^2) + \frac{C}{n \gamma^2} .
\end{eqnarray*}
where $C$ denotes a positive constant independent of $\epsilon$ and $n$.
\end{lemma}
PROOF. In a first time, remark that
\begin{eqnarray}
\E \sum_{|k|\leq N^{\star}} |\theta_k|^{2} \left( \left| \frac{\tilde \gamma_k}{\gamma_k} \right|^2 -1 \right)
& = & \E \sum_{|k|\leq N^{\star}} |\theta_k|^{2} |\gamma_k|^{-2} ( |\tilde \gamma_k - \gamma_k + \gamma_k|^2 -|\gamma_k|^{2}), \nonumber\\
& = & \E \sum_{|k|\leq N^{\star}} |\theta_k|^{2} |\gamma_k|^{-2} \left\lbrace |\tilde \gamma_k - \gamma_k|^2 +2 Re((\tilde\gamma_k -\gamma_k)\bar\gamma_k) \right\rbrace.
\label{eq:E11}
\end{eqnarray}
Let $N\in \mathbb{N}$ be a deterministic bandwidth. Since $\E \tilde \gamma_k = \gamma_k$ for all $k\in\mathbb{N}$, we can write that
\begin{eqnarray*}
\lefteqn{\E \sum_{|k|\leq N^{\star}} |\theta_k|^{2} |\gamma_k|^{-2} Re((\tilde\gamma_k -\gamma_k)\bar\gamma_k)}\\
& = & \E \sum_{|k|\in \lbrace N \dots N^{\star} \rbrace} |\theta_k|^{2} |\gamma_k|^{-2} Re((\tilde\gamma_k -\gamma_k)\bar\gamma_k) ,\\
& \leq & \E \left| \sum_{|k|\in \lbrace N \dots N^{\star} \rbrace} |\theta_k|^{2} |\gamma_k|^{-2} Re((\tilde\gamma_k -\gamma_k)\bar\gamma_k) \right| ,\\
& \leq & \E \sum_{k\in \mathbb{Z}} \left| (\1_{\lbrace |k| \leq N^{\star} \rbrace} - \1_{\lbrace |k| \leq N \rbrace})    |\theta_k|^{2} |\gamma_k|^{-2} Re((\tilde\gamma_k -\gamma_k)\bar\gamma_k) \right|.
\end{eqnarray*}
Using simple algebra
\begin{eqnarray*}
\left| \1_{\lbrace |k| \leq N^{\star} \rbrace} - \1_{\lbrace |k| \leq N \rbrace} \right|
& = & \left| (\1_{\lbrace |k| \leq N^{\star} \rbrace} + \1_{\lbrace |k| \leq N \rbrace})(\1_{\lbrace |k|\leq N^{\star} \rbrace} - \1_{\lbrace |k| \leq N \rbrace}) \right|,\\
& = &  (\1_{\lbrace |k| \leq N^{\star} \rbrace} + \1_{\lbrace |k| \leq N \rbrace}) \left| \1_{\lbrace |k| > N^{\star} \rbrace} - \1_{\lbrace |k| > N \rbrace} \right|,\\
& \leq &  \1_{\lbrace |k| > N^{\star} \rbrace}\1_{\lbrace |k| \leq N \rbrace} + \1_{\lbrace |k| \leq N^{\star} \rbrace}\1_{\lbrace |k| > N \rbrace}.
\end{eqnarray*}
For all $\gamma>0$, using the Cauchy-Schwartz and Young inequalities, we obtain
\begin{eqnarray}
\lefteqn{\E \sum_{|k|\leq N^{\star}} |\theta_k|^{2} |\gamma_k|^{-2} Re((\tilde\gamma_k -\gamma_k)\bar\gamma_k)} \nonumber \\
& \leq & \E \sum_{k\in \mathbb{Z}} \1_{\lbrace |k| > N^{\star} \rbrace}\1_{\lbrace |k| \leq N \rbrace} |\theta_k|^{2} |\gamma_k|^{-2} Re((\tilde\gamma_k -\gamma_k)\bar\gamma_k) \nonumber \\
& & \hspace{2cm} + \E \sum_{k\in \mathbb{Z}} \1_{\lbrace |k| \leq N^{\star} \rbrace}\1_{\lbrace |k| > N \rbrace} |\theta_k|^{2} |\gamma_k|^{-2} Re((\tilde\gamma_k -\gamma_k)\bar\gamma_k) \nonumber \\
& \leq & \gamma \E \sum_{|k|> N^{\star}} |\theta_k|^2 + \gamma \sum_{|k|> N} |\theta_k|^2 + \gamma^{-1} \E \sum_{|k|\leq N} |\theta_k|^2 |\gamma_k|^{-2} |\tilde\gamma_k - \gamma_k|^2 \nonumber \\
& & \hspace{3cm} + \gamma^{-1} \E \sum_{|k|\leq N^{\star}} |\theta_k|^2 |\gamma_k|^{-2} |\tilde\gamma_k - \gamma_k|^2.
\label{eq:E12}
\end{eqnarray}
Hence, from (\ref{eq:E11}) and (\ref{eq:E12})
\begin{eqnarray*}
\E \sum_{|k|\leq N^{\star}} |\theta_k|^{2} \left( \left| \frac{\tilde \gamma_k}{\gamma_k} \right|^2 -1 \right)
& \leq & (1+\gamma^{-1}) \E \sum_{|k|\leq N^{\star}} |\theta_k|^2 |\gamma_k|^{-2} |\tilde\gamma_k - \gamma_k|^2 + \gamma \E \sum_{|k|> N^{\star}} |\theta_k|^2 \\
& & \hspace{1.5cm} + \gamma  \sum_{|k|> N} |\theta_k|^2 + \gamma^{-1} \E \sum_{|k|\leq N} |\theta_k|^2 |\gamma_k|^{-2} |\tilde\gamma_k - \gamma_k|^2.
\end{eqnarray*}
A direct application of Lemma \ref{lemme:1} provides, for all $K>0$
\begin{eqnarray*}
\E \sum_{|k|\leq N^{\star}} |\theta_k|^{2} \left( \left| \frac{\tilde \gamma_k}{\gamma_k} \right|^2 -1 \right)
& \leq &  (1+\gamma^{-1})K \frac{\log^2(n)}{n} \E \sum_{|k|\leq N^{\star}} |\theta_k|^2 |\gamma_k|^{-2} + \gamma \E \sum_{|k|> N^{\star}} |\theta_k|^2 \\
& & \hspace{0.5cm} + \gamma  \sum_{|k|> N} |\theta_k|^2 + \frac{\gamma^{-1}}{n} \sum_{|k|\leq N} |\theta_k|^2 |\gamma_k|^{-2}(1-|\gamma_k|^2) + \frac{C}{nK}.
\end{eqnarray*}
Just set $K=\gamma^2$ in order to conclude the proof of Lemma \ref{lemme:2}.
\begin{flushright}
$\Box$
\end{flushright}

\begin{lemma}
\label{lemme:3}
Let $N^{\star}$ the bandwidth defined in (\ref{def:fenetre}). For all deterministic bandwidth $N$ and $0<\gamma<1$, we have
\begin{eqnarray*}
& & \frac{2\epsilon}{\sqrt{n}} \E \sum_{|k|\leq N^{\star}} |\gamma_k|^{-2} Re(\theta_k \tilde\gamma_k \bar \xi_k)
\leq 3\gamma \left\lbrace \sum_{|k|>N_0} |\theta_k|^2 + \frac{\epsilon^2}{n} \sum_{|k|\leq N_0 } |\gamma_k|^{-2} \right\rbrace \\
& & \hspace{2cm}+ 3\gamma \log^2 (n)\E \left\lbrace \sum_{|k|>N^{\star}} |\theta_k|^2 + \frac{\epsilon^2}{n} \sum_{|k|\leq N^{\star} } |\gamma_k|^{-2} \right\rbrace + \frac{C}{n} + \frac{C\epsilon^2}{\gamma^{4\beta+1}}.
\end{eqnarray*}
\end{lemma}
PROOF. In the following, we will use the inequality:
$$ P \left( \bigcup_{k=1}^{m_0} \left\lbrace \frac{1}{2} \leq \left| \frac{\tilde\gamma_k}{\gamma_k} \right| \leq 2 \right\rbrace \right) \leq \exp(-\log^{1+\tau} n),$$
for some $\tau>0$, wich can be proved using a Bernstein type inequality. Then, for all $\gamma>0$, using the above result and inequality (4.31) of \cite{riskhull}, we obtain
\begin{eqnarray*}
\frac{2\epsilon}{\sqrt{n}} \E \sum_{|k|\leq N^{\star}} |\gamma_k|^{-2} Re(\theta_k \tilde\gamma_k \bar \xi_k)
& \leq & \gamma \left\lbrace \sum_{|k|>N_0} |\theta_k|^2 + \frac{\epsilon^2}{n} \E \sum_{|k| \leq N_0} |\gamma_k|^{-4} |\tilde\gamma_k |^2 \right\rbrace \\
&  & + \gamma \E \left\lbrace \sum_{|k|>N^{\star}} |\theta_k|^2 + \frac{\epsilon^2}{n} \sum_{|k| \leq N^{\star}} |\gamma_k|^{-4} |\tilde\gamma_k |^2 \right\rbrace + \frac{C\epsilon^2}{n}\frac{1}{\gamma^{4\beta+1}}.
\end{eqnarray*}
In order to prove the above inequality, we use the inequality (4.31) of
\cite{riskhull} and 
Since $\E \tilde \gamma_k = \gamma_k$,
\begin{eqnarray*}
\frac{\epsilon^2}{n} \E \sum_{|k| \leq N} |\gamma_k|^{-4} |\tilde\gamma_k |^2
& = & \frac{\epsilon^2}{n} \sum_{|k| \leq N} |\gamma_k|^{-4} \left\lbrace \E |\tilde\gamma_k |^2 - |\gamma_k|^2 + |\gamma_k|^2 \right\rbrace,\\
& = & \frac{\epsilon^2}{n} \sum_{|k| \leq N} |\gamma_k|^{-4} \left\lbrace |\gamma_k |^2 + \mathrm{Var}(\tilde\gamma_k) \right\rbrace,\\
& = &  \frac{\epsilon^2}{n} \sum_{|k|\leq N} |\gamma_k|^{-4} \left\lbrace |\gamma_k|^2 + \frac{1}{n} \right\rbrace \leq 2 \frac{\epsilon^2}{n} \sum_{|k|\leq N_0} |\gamma_k|^{-2}.
\end{eqnarray*}
The same kind of inequality can be obtained with the random bandwidth $N^{\star}$. Indeed,
\begin{eqnarray*}
\frac{\epsilon^2}{n} \E \sum_{|k|\leq N^{\star}} |\gamma_k|^{-4}|\tilde\gamma_k|^2
& = & \frac{\epsilon^2}{n} \E \sum_{|k|\leq N^{\star}} |\gamma_k|^{-4}|\tilde\gamma_k - \gamma_k + \gamma_k|^2,\\
& \leq & \frac{\epsilon^2}{n} \E \sum_{|k|\leq N^{\star}} |\gamma_k|^{-4} \left\lbrace 2|\tilde\gamma_k - \gamma_k|^2 + 2|\gamma_k|^2 \right\rbrace,\\
& \leq & \frac{2\epsilon^2}{n} \E \sum_{|k|\leq N^{\star}} |\gamma_k|^{-2} + \frac{2\epsilon^2}{n} \E \sum_{|k|\leq N^{\star}} |\gamma_k|^{-4} |\tilde\gamma_k - \gamma_k|^2.
\end{eqnarray*}
Using the same algebra as in the proof of Lemma 5.1, we obtain, for all $Q>0$:
\begin{eqnarray*}
\lefteqn{\frac{\epsilon^2}{n} \E \sum_{|k|\leq N^{\star}} |\gamma_k^{-4}| |\tilde\gamma_k - \gamma_k|^2}\\
& = & Q \frac{\epsilon^2}{n} \E \sum_{|k|\leq N^{\star}} |\gamma_k|^{-4} + \frac{\epsilon^2}{n} \E \sum_{|k|\leq N^{\star}} |\gamma_k|^{-4} \left\lbrace |\tilde\gamma_k - \gamma_k|^2 - Q \right\rbrace \1_{\lbrace |\tilde\gamma_k - \gamma_k|^2>Q \rbrace},\\
& \leq & Q \frac{\epsilon^2}{n} \E \sum_{|k|\leq N^{\star}} |\gamma_k|^{-4} + \frac{C\epsilon^2 n}{\log^2(n)} \E \sum_{|k|\leq N^{\star}} |\tilde\gamma_k - \gamma_k|^2 \1_{\lbrace |\tilde\gamma_k - \gamma_k|^2>Q \rbrace},\\
& \leq & Q \frac{\epsilon^2}{n} \E \sum_{|k|\leq N^{\star}} |\gamma_k|^{-4} + \frac{C\epsilon^2}{\log^4(n)}  e^{-Qn/4}.
\end{eqnarray*}
Setting $Q=n^{-1}\log^2(n)$, we obtain,
\begin{eqnarray*}
\frac{\epsilon^2}{n} \E \sum_{|k|\leq N^{\star}} |\gamma_k^{-4}| |\tilde\gamma_k - \gamma_k|^2
& \leq & \frac{\epsilon^2}{n} \E \sum_{|k|\leq N^{\star}} |\gamma_k|^{-2} \frac{|\gamma_k|^{-2} \log^2(n)}{n} + \frac{C\epsilon^2}{n},\\
& \leq & \frac{\epsilon^2}{n} \E \sum_{|k|\leq N^{\star}} |\gamma_k|^{-2} + \frac{C\epsilon^2}{n}.
\end{eqnarray*}
This concludes the proof.
\begin{flushright}
$\Box$
\end{flushright}

\bibliography{rand_shift_oracle}
\bibliographystyle{alpha}

\end{document}